\input amstex
\documentstyle{amsppt}
\magnification=\magstephalf
\loadbold
\NoRunningHeads
\NoBlackBoxes
\TagsOnRight
\topmatter
\title
\nofrills
Invariant Hyperk\"ahler Structures on the
Cotangent Bundles of Hermitian Symmetric Spaces
\endtitle
\author
I.V.~Mykytyuk
\endauthor
\address
Institute of Applied Mathematics and Fundamental Sciences,
National University ``L'viv Polytechnica'',
S. Bandery Str., 12,
79013 L'viv, Ukraine
\endaddress
\email
viva\@iapmm.lviv.ua
\endemail
\thanks
This research was partially supported by the
German National Science Foundation
(Sonderforschungsbereich 237).
\endthanks
\keywords
hyper-K\"ahlerian structure, Hermitian symmetric space
\endkeywords
\subjclass
32Q15, 37J15
\endsubjclass
\abstract
Let
$G/K$ be an irreducible Hermitian symmetric space of compact
type with the standard homogeneous complex
structure. Then the real symplectic manifold
$(T^*(G/K),\Omega)$ has the natural complex structure
$J^-$. We construct all
$G$-invariant K\"ahler structures
$(J,\Omega)$ on homogeneous domains in
$T^*(G/K)$ anticommuting with
$J^-$. Each such a hypercomplex structure, together with a
suitable metric, defines a hyper-K\"ahlerian structure.
As an application, we obtain a new proof of
the Harish-Chandra and Moore theorem
for Hermitian symmetric spaces.

Bibliography: 13 titles.
\endabstract
\endtopmatter

\document

\head
{\bf \S 1. Introduction}
\endhead

Let
$M=G/K$ be an irreducible Hermitian symmetric
space of compact type with a homogeneous metric
${\bold g}_{\text M}$. Since
$M$ is a homogeneous complex manifold, its cotangent bundle
$T^*M$ has a natural complex structure. Using
${\bold g}_{\text M}$ we can identify the cotangent and tangent
bundles and thus obtain a complex structure on
$TM$, with respect to which the zero section
$M\subset TM$ is complex. This structure
$J^-$ is different from the standard complex structure
$J^+$ on $TM$ induced by that on
$M$.

On the other hand, the cotangent bundle
$T^*M\simeq TM$ is a symplectic manifold
with the canonical symplectic form
$\Omega$. In this paper we make an explicit
description of all $G$-invariant K\"ahler structures
$(J,\Omega)$ (with the K\"ahler form
$\Omega$) on homogeneous domains
$D\subset TM$  anticommuting with
$J^-$ (Theorem~4.12). In fact, each resulting hypercomplex
structure, together with the suitable metric
${\bold g}$, defines a hyper-K\"ahlerian structure.

If the domain
$D$ contains the zero section
$M$, the restriction of the hyper-K\"ahlerian metric
${\bold g}$ to
$M$ is the given homogeneous metric
${\bold g}_{\text M}$ up to a constant
multiplier (one makes this multiplier
$=1$ using for the identification of
$T^*M$ and $TM$ a homogeneous metric on
$M$ proportional to
${\bold g}_{\text M}$). Such hyper-K\"ahlerian metrics have been
constructed in~\cite{Bu} using twistor methods and case by case
the classification of symmetric spaces, in~\cite{Bi} using Nahm's
equations and in~\cite{DSz} (for spaces of classical groups) using
deformation of the so-called adapted complex structure on
$TM$. In~\cite{BG1} Biquard and Gauduchon found explicit formulas
for these hyper-K\"ahlerian metrics
in terms of some operator-functions
$P:{\frak m}\to\text{End}({\frak m})$ on the space
${\frak m}\simeq T_o(G/K)$, where
$o=\{K\}$. These hyper-K\"ahlerian structures are global ones. Our
additional structures are not defined on the zero section
$M$. So we cannot talk about a restriction
of the corresponding hyper-K\"ahlerian metric to
$M$ as in~\cite{BG1}. Nevertheless, our expressions for
$P$ and potential functions generalize the corresponding formulas
of~\cite{BG1,BG2}.

For proofs in~\cite{DSz,BG1,BG2} they used the decomposition of
$T(TM)$ between horizontal and vertical
directions, induced by the Levi-Civita connection of
$M$. Our approch is based on the fact that
$T(G/K)$ is a reduced manifold for
the (right) Hamiltonian action of $K$ on
$TG$. We can substantially simplify matters by working as
in~\cite{My1,My2} in the trivial vector bundle
$G\times{\frak m}$ which is a level surface for the corresponding
moment map. So we use the natural homogeneous decomposition
of $T(G\times{\frak m})$ usual for the Lie algebras theory. As an
application we obtain a new simple proof of the well-known
Harish-Chandra and Moore theorem about restricted root systems of
Hermitian symmetric spaces.

The part of this work was done while the author was visiting
the Ruhr--University (Bochum, Germany) in November--December,
2001. The author would like to thank Prof. A. Huckleberry for
support and hospitality. Besides, I would like to express my
gratitude to Prof. A.M. Stepin for very helpful discussions.

\head
{\bf\S 2. $\bold G$-invariant K\"ahler structures on $\bold {T(G/K)}$}
\endhead

{\bf 2.1. Anticommuting structures.}
We recall some facts on hypercomplex and hyper-K\"ahlerian
structures (see for example~\cite{BG1,Ob,Hi}). Let
$N$ be a smooth real manifold with a complex structure
$J$ and a symplectic 2-form
$\omega$ (all objects in this paper are smooth,
unless otherwise indicated).
For any vector bundle $L$ on
$N$ denote by
${\Gamma}L$ the set of its smooth sections. For the tensor
$J$ denote by
$F(J)\subset T^{\Bbb C} N$ its (involutive) complex subbundle of
$(0,1)$-vectors, i.e.
${\Gamma}F(J)=\{X+iJX,X\in {\Gamma}(TN)\}$. We need some
definitions.

\definition{Definition 2.1}
The pair
$(J,\omega)$ is a K\"ahler structure on
$N$ if
\roster
\item the closed 2-form
$\omega$ is invariant with respect to $J$:
$\omega(JX,JY)=\omega (X,Y)$,
$\forall X,Y\in{\Gamma}(TM)$;
\item the bilinear form
${\bold g}={\bold g}(J,\omega)$, where
${\bold g}(X,Y)\overset\text{def}\to{=}\omega(JX,Y)$,
is symmetric and positive definite.
\endroster
We will denote
such a K\"ahler structure also by the pair
$(J,{\bold g})$ because
$\omega(X,Y)={\bold g}(-JX,Y)$, i.e.
$\omega=\omega(J,{\bold g})$.
\enddefinition

Let
$\Pi: \tilde N\to N$ be a submersion of a manifold
$\tilde N$ onto $N$ and
${\Cal K}\subset T\tilde N$ be the kernel of
$\Pi_*$. Let
${\Cal T}$ be some complementary subbundle to
${\Cal K}$ in $T\tilde N$, i.e.
${\Cal K}\oplus{\Cal T}=T\tilde N$. For the complex structure
$J$ on $N$ there exists a unique (smooth)
$(1,1)$-tensor $\tilde J$ on
$\tilde N$ such that
$$
\tilde J({\Cal T})={\Cal T},
\qquad \tilde J({\Cal K})=0,
\qquad
\Pi_*\circ\tilde J=J\circ\Pi_*
\quad\text{on}\quad {\Cal T}.
\tag 2.1
$$
Similarly, for the $(0,1)$-subbundle
$F(J)$ there exists a unique complex subbundle
${\Cal F}(J)\subset T^{\Bbb C}\tilde N$ containing the kernel
${\Cal K}$ and such that
$\Pi_*({\Cal F}(J))=F(J)$. It is clear that this subbundle is
involutive.

\proclaim{Lemma 2.2}
The form
$\omega$ is invariant with respect to
$J$ iff the (1,1)-tensor
$J$ is skew-symmetric with respect to the form
$\omega$: $\omega(JX,Y)=\omega (X,-JY)$,
$\forall X,Y\in{\Gamma}(TM)$. The tensor
$J$ is skew-symmetric (with respect to
$\omega$) iff
$\omega(F(J),F(J))=0$, and symmetric iff
$\omega(F(J),\overline{F(J)})=0$.
\endproclaim
\demo{Proof}
Taking into account that $J^2=-1$ and
$$
\omega(X+iJX,Y\pm iJY)=
i\Bigl[\omega(JX,Y)\pm \omega(X,JY)\Bigr]+
\Bigl[\omega(X,Y)\mp \omega(JX,JY)\Bigr]
$$
we obtain the assertions of the lemma.
{}\hfill$\square$
\enddemo

Observe the following fact:
\proclaim{Corollary 2.2.1}
The tensor
$J$ is skew-symmetric with respect to the form
$\omega$ iff
$(\Pi^*\omega)({\Cal F}(J),{\Cal F}(J))=0$, and symmetric iff
$(\Pi^*\omega)({\Cal F}(J),\overline{{\Cal F}(J)})=0$.
\endproclaim

The following assertion is well-known (see~\cite{GS, Lemma~4.3}).
\proclaim{Lemma 2.3}
The pair
$(J,\omega)$ is a K\"ahler structure on
$N$ iff $\omega(F(J),F(J))=0$ and
$-i\omega(Z,\overline{Z})>0$ for all non-zero local vector-fields
$Z\in {\Gamma}_{loc}F(J)$.
\endproclaim

\definition{Definition 2.4}
A pair
$(J_1,J_2)$ formed by two anticommuting complex structures
$J_1$ and
$J_2$ is a hypercomplex structure on
$N$. Then
$J_3=J_1J_2$ is also a complex structure on
$N$ (for a proof see~\cite{Ob}).
\enddefinition

\remark{Remark 2.5}
Almost-complex tensors $J_1$ and $J_2$ on $N$ are anticommuting
iff so are the tensors $\tilde J_1$ and $\tilde J_2$ on $\tilde N$.
\endremark

\definition{Definition 2.6}~\cite{BG1}
A triple
$({\bold g}, J_1,J_2)$ formed by a Riemannian metric
$\bold g$ and two anticommuting complex structures
$J_1$ and $J_2$ on
$N$ is a hyper-K\"ahlerian structure on
$N$ whenever the pairs $(J_1,{\bold g})$ and
$(J_2,{\bold g})$ are K\"ahler. Then the pair
$(J_3=J_1J_2,{\bold g})$ is a K\"ahler one as well.
\enddefinition

For a hyper-K\"ahlerian structure
$({\bold g}, J_1,J_2)$, let us denote by
$\omega_j$ the K\"ahler form corresponding to
$J_j$,
$j=1,2,3$. It is clear that this hyper-K\"ahlerian
structure is determined by any pair
$(J_k,(J_j,\omega_j))$,
$k\ne j$. Since $J_k$ and
$J_j$ anticommute, the tensor
$J_k$ is symmetric with respect to the form
$\omega_j$ (while
$J_j$ is skew-symmetric). It is clear that
$\omega_k(X,Y)=\omega_j(J_kJ_j X,Y)$,
$X,Y\in {\Gamma}(TN)$
($k\ne j$).
The following simple lemma defines a hyper-K\"ahlerian structure in
these terms.
\proclaim{Lemma 2.7}
Let $\bigl(J',(J,\omega)\bigr)$ be a pair, where
$J',J$ are anticommuting complex structures and
$(J,\omega)$ is a K\"ahler structure on
$N$ with the corresponding Hermitian metric
${\bold g}={\bold g}(J,\omega)$. Then the triple
$({\bold g}, J',J)$ is a hyper-K\"ahlerian structure on
$N$ iff
\roster
\item the tensor
$J'$ is symmetric with respect to
$\omega$;
\item the 2-form $\omega'$,
where $\omega'(X,Y)\overset\text{def}\to{=} \omega(J'JX,Y)$
$\forall X,Y\in {\Gamma}(TN)$, is closed.
\endroster
\endproclaim
\demo{Proof}
Since $J'J=-JJ'$, from (1) it follows that
$$
\omega'(X,Y)\overset\text{def}\to{=}
\omega(J'JX,Y)=
\omega(X,-JJ'Y)=
\omega(JJ'Y,X)=
-\omega'(Y,X),
$$
i.e. the bilinear form $\omega'$ is skew-symmetric.
Now to prove the lemma it is sufficient to verify that
the form $\omega'$ is $J'$-invariant
and ${\bold g}(X,Y)=\omega'(J'X,Y)$.
These properties follow from the following chains of equations:
$$
\omega'(J'X,J'Y)\overset\text{def}\to{=}
\omega(J'JJ'X,J'Y)=
\omega(JX,J'Y)=
\omega(J'JX,Y)\overset\text{def}\to{=}
\omega'(X,Y)
$$
and
$$
{\bold g}(X,Y)\overset\text{def}\to{=}
\omega(JX,Y)=
\omega(J'JJ'X,Y)\overset\text{def}\to{=}
\omega'(J'X,Y).
$$
\enddemo

Since the kernels of the $(1,1)$-tensors $\tilde J'$,
$\tilde J$ and the forms $\Pi^*\omega'$, $\Pi^*\omega$
coincide with ${\Cal K}$, it follows from the definition
of $\tilde J'$ and $\tilde J$ (see~(2.1))
that
$$
(\Pi^*\omega')(\tilde X,\tilde Y)=
(\Pi^*\omega)(\tilde J'\tilde J\tilde X,\tilde Y),
\qquad
\forall \tilde X,\tilde Y\in {\Gamma}(T\tilde N).
\tag 2.2$$

{\bf 2.2. $\bold G$-invariant K\"ahler structures
$\bold{(J(P),\boldsymbol\Omega)}$.}
Let $M=G/K$ be a symmetric space with
a real reductive connected Lie group
$G$ and a compact subgroup $K$. Let
${\frak g}$ and
${\frak k}$ be the Lie algebras of the groups
$G$ and
$K$ respectively,
$$
{\frak g}={\frak k}\oplus{\frak m},
\qquad
[{\frak k},{\frak m}]\subset{\frak m}
\qquad
[{\frak m},{\frak m}]\subset{\frak k}.
\tag 2.3$$
Suppose that there is a nondegenerate
$\operatorname{Ad} G$-invariant bilinear form
$\langle ,\rangle$ on
${\frak g}$ such that its restriction
$\langle ,\rangle|{\frak m}$ is a positive definite form
and ${\frak k}\bot{\frak m}$.
This form defines
$G$-invariant Riemannian metric
${\bold g}_{\text M}$ on
$M=G/K$. The metric
${\bold g}_{\text M}$ identifies the cotangent bundle
$T^*M$ and the tangent bundle
$TM$ and thus we can also talk about the canonical 1-form
$\theta$ on $TM$. The form
$\theta$ and the symplectic form
$\Omega\overset\text{def}\to{=} d\theta$ are
$G$-invariant with respect to the natural action of
$G$ on
$TM$.

Since
${\frak g}={\frak k}\oplus{\frak m}$ is
$\operatorname{Ad}(K)$-invariant (orthogonal) splitting of
${\frak g}$, we can consider a trivial vector bundle
$G\times {\frak m}$ with the two Lie group
actions (which commute) on it: the left
$G$-action,
$l_h:(g,w)\mapsto (hg,w)$ and the right
$K$-action
$r_k:(g,w)\mapsto (gk,\operatorname{Ad}_{k^{-1}}w)$. Let
$\pi: G\times {\frak m}\to G\times_K {\frak m}$ be the
natural projection. It is well known that
$G\times_K {\frak m}$ and
$TM$ are isomorphic. Using the corresponding
$G$-equivariant diffeomorphism
$\varphi: G\times_K {\frak m}\to TM$,
$[(g,w)]\mapsto \frac{d}{dt}\bigr|_0 g\exp(tw)K$
and the projection $\pi$ define the
$G$-equivariant submersion
$\Pi: G\times {\frak m}\to TM$,
$\Pi=\varphi\circ\pi$.
Let
$\xi^l$ be the left-invariant vector field on the Lie group
$G$ defined by a vector
$\xi\in{\frak g}$. By~\cite{My1, Lemma~2.3}
$$
(\Pi^*\theta)_{(g,w)}\bigl(\xi^l(g), u\bigr)=
\langle w,\xi\rangle,\tag 2.4
$$
$$
(\Pi^*\Omega)_{(g,w)}\bigl((\xi_1^l(g), u_1),
(\xi_2^l(g), u_2)\bigr)=
\langle\xi_2,u_1\rangle
-\langle\xi_1,u_2\rangle
-\langle w, [\xi_1,\xi_2]\rangle, \tag 2.5
$$
where
$g\in G, w\in {\frak m}$,
$\xi,\xi_1,\xi_2 \in {\frak g}$,
$u,u_1,u_2\in {\frak m}=T_w {\frak m}$.
Since $\Omega$ is a symplectic form,
the kernel ${\Cal K}\subset T(G\times {\frak m})$ of the 2-form
$\Pi^*\Omega$ is the kernel of
$\Pi_*$.

Let $D$ be an open connected
$G$-invariant subset of $TM$. Denote by
$W$ a unique
$\operatorname{Ad}(K)$-invariant open subset of
${\frak m}$ such that
$\Pi^{-1}(D)=G\times W$. Let
$\text{Eqv}(W)$ be the set of all smooth
$K$-equivariant mappings
$A: W\to \text{End}({\frak m}^{\Bbb C})$,
$w\mapsto A_w$, i.e. for which
$$
\operatorname{Ad}_k\circ A_w\circ
\operatorname{Ad}_{k^{-1}}= A_{\operatorname{Ad}_k w}
\quad\text{on}\quad {\frak m}\quad
\text{for all}\quad w\in W,\ k\in K.
\tag 2.6$$
Denote by $\text{Alm}(W)$ the set of all
$P\in\text{Eqv}(W)$ such that the operator
$P_w:{\frak m}^{\Bbb C}\to{\frak m}^{\Bbb C}$ and its real part
$\operatorname{Re} P_w:{\frak m}\to{\frak m}$
are nondegenerate for each
$w\in W$. Such a $K$-equivariant mapping
$P\in\text{Alm}(W)$ determines a complex (left)
$G$-invariant subbundle
${\Cal F}(P)\subset T^{\Bbb C}(G\times W)$ generated by (left)
$G$-invariant vector fields $\xi^L$,
$\xi\in{\frak m}$ and $\zeta^L$,
$\zeta\in{\frak k}$ on
$G\times W$, where
$$
\xi^L(g,w)=\bigl(\xi^l(g),iP_w(\xi)\bigr),
\qquad
\zeta^L(g,w)=\bigl(\zeta^l(g),[w,\zeta]\bigr).
$$
The subbundle ${\Cal F}(P)$ is (right)
$K$-invariant by~(2.6) and because the vector fields
$\zeta^L$, $\zeta\in {\frak k}$ span the (right)
$K$-invariant subbundle (kernel) ${\Cal K}$. Therefore
$F(P)\overset\text{def}\to{=}\Pi_*({\Cal F(P)})$ is a
well-defined (smooth) complex subbundle of
$T^{\Bbb C} D$
(${\Cal K}^{\Bbb C}\subset{\Cal F}(P)$) such that
$F(P)+\overline{F(P)}=T^{\Bbb C}D$ ¨
$F(P)\cap\overline{F(P)}=0$. In other words, the mapping
$P$ determines an almost-complex structure
$J(P)$ on $D\subset TM$ with
$F(P)$ as the subbundle of its
$(0,1)$-vectors.

For a vector field $X\in{\Gamma}(TW)$ and
$A\in\text{Eqv}(W)$ denote by
$L_X A$ the derivative of $A$ along
$X$, i.e.
$(L_X A)_w\overset\text{def}\to{=}
(d/dt)_0 A_{w+tX(w)}\in\text{End}({\frak m}^{\Bbb C})$.
We extend this definition on complex
vector fields using linearity.
Each vector $\xi\in{\frak m}$ defines the vector field $P\xi$
on $W$ by $(P\xi)_w=P_w(\xi)$.
Now we want to present a result which will be effectively used
in a remaining part of the paper.

\proclaim{Proposition 2.8}~\cite{My1}
Let $M=G/K$ be a Riemannian symmetric space.
The almost-complex structure
$J(P)$ on the domain $D=\Pi(G\times W)\subset TM$ is
\roster
\item integrable iff
for all (fixed) vectors $\xi,\eta\in{\frak m}$ and $w\in W$
$$
(L_{P\xi} P)_w(\eta)-(L_{P\eta} P)_w(\xi)=
-{\bigl[w,[\xi,\eta]\bigr]};
\tag 2.7
$$
\item a K\"ahler structure with the K\"ahler form
$\Omega$ iff {\rm(1)} holds and for each
$w\in W$ the endomorphism
$P_w$ is symmetric with positive-definite real part
$\operatorname{Re} P_w$ {\rm(}with respect to the bilinear form
$\langle,\rangle$ on
${\frak m}${\rm)}.
\endroster

For any $G$-invariant K\"ahler structure $J$ on
$D$ with $\Omega$ as a K\"ahler form
there exists a unique
mapping $P\in\text{Alm}(W)$ such that
$J=J(P)$.
\endproclaim

Observe the following fact:
\proclaim{Corollary 2.8.1}
If the structure $J(P)$ is integrable then
so are $J(-P)$ and $J(\overline P)$.
\endproclaim

\proclaim{Corollary 2.8.2}
Suppose that
$(J(P),\Omega)$ is a K\"ahler structure on
$D$. If $0\in W$ and the Lie algebra
${\frak g}$ is simple then
$P_0=\psi_0\cdot\operatorname{Id}_{\frak m}$, where
$\operatorname{Re}\psi_0\in{\Bbb R}^+$.
\endproclaim
\demo{Proof} If
${\frak g}$ is a simple algebra,
${\frak m}$ is a simple
$\operatorname{Ad}(K)$-module~\cite{GG, (8.5.1)}. Since
$P_0$ is a symmetric endomorphism
which commutes with all endomorphisms
$\operatorname{Ad}_k|{\frak m}$,
$k\in K$ (condition~(2.6)),
$P_0=\psi_0\cdot\operatorname{Id}_{\frak m}$ for some
$\psi_0\in{\Bbb C}$.
{}\hfill$\square$
\enddemo

\head
{\bf\S 3. Invariant K\"ahler structures
on Hermitian symmetric spaces}
\endhead

We continue with the previous notations but
in this section it is assumed in addition that
$G/K$ is an Hermitian symmetric
space, i.e. there exists an endomorphism
$I:{\frak m}\to{\frak m}$ such that

(1)~$I^2=-\operatorname{Id}_{\frak m}$;

(2)~$\operatorname{Ad}_kI=I\operatorname{Ad}_k$ on
${\frak m}$, $\forall k\in K$;

(3)~the form $\langle ,\rangle|{\frak m}$ is
$I$-invariant; and
$$
[I\xi,I\eta]=[\xi,\eta], \qquad I[\zeta,\eta]=[\zeta,I\eta]
\qquad
\text{for all } \xi,\eta\in{\frak m},\ \zeta\in{\frak k}.
\tag 3.1
$$
Such a triple
$({\frak g},{\frak k},I)$ we will call an Hermitian orthogonal
symmetric Lie algebra. It follows from~(3.1) that
$[\xi^l+i(I\xi)^l,\eta^l+i(I\eta)^l]=0$,
$\forall \xi,\eta\in{\frak m}$. Thus the complex subbundle of
$T^{\Bbb C} G$, generated by vector-fields
$\xi^l+i(I\xi)^l$,
$\xi\in{\frak m}$, is involutive. Since this subbundle is left
$G$-invariant and right
$K$-invariant, its image under the canonical projection
$G\to G/K$ defines a
$G$-invariant complex structure on
$M=G/K$.

{\bf 3.1. Hypercomplex structures on the tangent
bundles of Hermitian symmetric spaces.}
Here we apply the general results of the previous
Section 2 in the special situation here.
The $G$-invariant complex structure on
$M$ induces the $G$-invariant
complex structures $J^+$ and $J^-$ on
$TM$ and
$T^*M\simeq TM$ respectively. It is clear that the subbundle
$F^\pm\subset T^{\Bbb C} (TM)$ of
$(0,1)$ vectors of
$J^\pm$ coincide with the subbundle
$\Pi_*({\Cal F}^\pm)$, where
${\Cal F}^\pm$ is the (left)
$G$-invariant and (right)
$K$-invariant subbundle of
$T^{\Bbb C}(G\times{\frak m})$:
$$
{\Cal F}^\pm(g,w)=\bigl\{\bigl(\xi^l(g)
+i(I\xi)^l(g),u\pm i(Iu)\bigr),\
\xi,u\in{\frak m}\bigr\}\oplus{\Cal K}(g,w).
\tag 3.2
$$

Fix some mapping
$P\in\text{Alm}(W)$. It is clear that the mapping
$PI$,
$(PI)_w\overset\text{def}\to{=} P_wI$
is also an element of the set
$\text{Alm}(W)$ because the group
$\operatorname{Ad}(K)|{\frak m}$ commutes elementwise with
$I$.
\proclaim{Lemma 3.1}
If $J(P)$,
$P\in\text{Alm}(W)$ is a complex structure then so is
$J(PI)$.
\endproclaim
\demo{Proof}
Since
$F(P)$ is an involutive subbundle and
$I$ is independent of
$w$, from Proposition~2.8 and relations~(3.1) it follows that
$$
(L_{P(I\xi)} P)_w(I\eta)-(L_{P(I\eta)} P)_w(I\xi)=
-\bigl[w,[I\xi,I\eta]\bigr]=-\bigl[w,[\xi,\eta]\bigr]
$$
on $W$ for any vectors
$\xi,\eta\in{\frak m}$, i.e. the subbundle
$F(PI)$ is also involutive.
{}\hfill$\square$
\enddemo

In order to describe the defined above complex structures in
terms of their almost-complex tensors, consider two (left)
$G$-invariant and (right)
$K$-invariant subbundles ${\Cal T}_h$ and
${\Cal T}_v$ of the tangent bundle
$T(G\times W)$ given by
$$
{\Cal T}_h(g,w)=\{(\xi^l(g),0),\ \xi\in{\frak m}\},
\quad
{\Cal T}_v(g,w)=\{(0,u),\ u\in{\frak m}=T_w W\}.
$$
Put
${\Cal T}={\Cal T}_h\oplus{\Cal T}_v$. Then
$T(G\times W)={\Cal K}\oplus{\Cal T}$. The mapping
$(\xi^l(g),0)\mapsto (0,\xi)$ determines
the canonical isomorphism of the spaces
${\Cal T}_h(g,w)$ and
${\Cal T}_v(g,w)$. Using this isomorphism, we obtain that the
$(1,1)$-tensors $\tilde J^\pm$ and
$\tilde J(P)$ on
$G\times W$ (see~(2.1)) at the point
$(g,w)$ are given by
$$
\tilde J^\pm_{(g,w)}|{\Cal T}=\left(
\matrix
I & 0\\
0 & \pm I
\endmatrix
\right),
\quad \tilde J_{(g,w)}(P)|{\Cal T}= \left(
\matrix
-R_w^{-1}S_w & -R_w^{-1}\\
R_w+S_wR_w^{-1}S_w & S_wR_w^{-1}
\endmatrix
\right),
\tag 3.3$$
where $R=\operatorname{Re} P$,
$S=\operatorname{Im} P$. Now it is easy to verify that the tensors
$\tilde J^-$ and
$\tilde J(P)$ are anticommuting iff
$RI=IR$ and
$SI=-IS$, i.e. by Remark~2.5
$$
J^-J(P)=-J(P)J^- \
\Longleftrightarrow \
RI=IR,\  SI=-IS \
\Longleftrightarrow \
IP=\overline PI.
\tag 3.4$$
Considering now the almost complex structures
$J(P)$ and $J(PI)$ with a real mapping
$P\in\text{Alm}(W)$ (i.e. with
$\operatorname{Im} P=0$), we obtain that
$$
\tilde J_{(g,w)}(P)|{\Cal T}= \left(
\matrix
0 & -P_w^{-1}\\
P_w & 0
\endmatrix
\right),
\qquad \tilde J_{(g,w)}(PI)|{\Cal T}= \left(
\matrix
0 & IP_w^{-1}\\
P_wI & 0
\endmatrix
\right).
$$
In this case
$J(P)J(PI)=-J(PI)J(P)$. We have established the following result.
\proclaim{Proposition 3.2}
Let $J(P)$,
$P\in\text{Alm}(W)$ be a complex structure on
$D$ such that $\operatorname{Im} P=0$. Then the pair
$\bigl(J(P),J(PI)\bigr)$ is a hypercomplex structure on
$D$.
\endproclaim

{\bf 3.2. Hyperk\"ahler structures on the tangent
bundles of Hermitian symmetric spaces.}
In this subsection we study properties of the pair
$\bigl(J^-,(J(P),\Omega)\bigr)$, where
$(J(P),\Omega)$ is a K\"ahler structure on the domain
$D\subset TM$.
\proclaim{Theorem 3.3}
Let
$(J(P),\Omega)$, $P\in\text{Alm}(W)$ be a K\"ahler structure on
$D$ with the Hermitian metric
${\bold g}={\bold g}(J(P),\Omega)$. Then the triple
$({\bold g},J^-, J(P))$ is a hyper-K\"ahlerian structure {\rm(}on
$D${\rm)} iff $IP=\overline PI$ on
$W$.
\endproclaim
\demo{Proof}
By~(3.4)
the complex structures
$J^-$ and $ J(P)$ anticommute iff
$IP=\overline PI$. For the pair
$(J^-, J(P))$ of anticommuting complex
structures the almost-complex structure
$J'=J^-J(P)$ is integrable~\cite{Ob}. Therefore the triples
$({\bold g},J^-, J(P))$ and
$({\bold g},J^-J(P),\allowbreak J(P))$ are hypercomplex structures
simultaneously. Thus to prove the theorem
it is sufficient to show that for the pair
$(J',(J,\Omega))$, where
$J=J(P)$, conditions (1) and (2) of Lemma~2.7 hold.

Since by definition the form $\Omega$ is
$J$-invariant and
$JJ'=-J'J$, we derive the identities
$\Omega(JJ'X,JY)=\Omega(J'X,Y)$ and
$\Omega(X,JJ'(JY))=\Omega(X,J'Y)$, where
$X,Y\in{\Gamma}(TD)$. In other words, the tensor
$J'$ is symmetric with respect to
$\Omega$ iff so is $JJ'$. But the tensor
$JJ'=J^-$ is symmetric with respect to
$\Omega$ because
$(\Pi^*\Omega)({\Cal F}^-,\overline{{\Cal F}^-})=0$ (see
Corollary~2.2.1). Indeed, using the relation
$[{\frak m},{\frak m}]\bot{\frak m}$, property (3) of
$I$ and definitions~(2.5),~(3.2) of
$\Pi^*\Omega$ and
${\Cal F}^-$, we obtain that for all
$\xi_1,\xi_2,u_1,u_2\in{\frak m}$
$$\align
(\Pi^*\Omega)
\bigl((\xi_1^l+i(I\xi_1)^l, u_1-i(Iu_1)),
(\xi_2^l-i(I\xi_2)^l, u_2+i(Iu_2))\bigr)&=\\
{}=\langle\xi_2-iI\xi_2,u_1-iIu_1\rangle
-\langle\xi_1+iI\xi_1,u_2+iIu_2\rangle&=0-0.
\endalign$$

Define the 1-form $\theta'$ and the tensor
$\Omega'$ on $D$ putting for
$X,Y\in {\Gamma}(TD)$
$$
\theta'(X)\overset\text{def}\to{=}\theta(-J^-X)
\quad\text{and}\quad
\Omega'(X,Y)\overset\text{def}\to{=}
\Omega(-J^-X,Y)=\Omega(J'JX,Y).
\tag 3.5$$
Since we already have proved that
$J^-$ is symmetric with respect to
$\Omega$, we have
$\Omega'(X,Y)=-\Omega'(Y,X)$. To prove that the form
$\Omega'$ is closed we will show that
$\Omega'=d\theta'$. By~(3.3) and by definition of the form
$\theta$
$$
(\Pi^*\theta')_{(g,w)}(\xi^l(g),u)
=\langle w,-I\xi\rangle
=\langle Iw,\xi\rangle,
\qquad
\xi\in{\frak g},\ u\in{\frak m}.
$$
Applying the well-known formula
$d\theta'(X,Y)=X\theta'(Y)-Y\theta'(X)-\theta'([X,Y])$ to the form
$\Pi^*\theta'$ we obtain that for
$\xi_1,\xi_2,u_1,u_2\in{\frak m}$
$$
d(\Pi^*\theta')_{(g,w)}
((\xi_1^l(g), u_1), (\xi_2^l(g), u_2))=
\langle\xi_2,Iu_1\rangle
-\langle\xi_1,Iu_2\rangle,
\tag 3.6$$
because $\langle Iw, [\xi_1,\xi_2]\rangle=0$. But using
expression~(3.3) for almost-complex tensor
$J^-$, we derive from~(2.2) the following formula for
$(\Pi^*\Omega')|({\Cal T}_h\oplus{\Cal T}_v)$
$$\align
&(\Pi^*\Omega')((\xi_1^l, u_1), (\xi_2^l, u_2))=
(\Pi^*\Omega)\Bigl(((-I\xi_1)^l, Iu_1), (\xi_2^l, u_2)\Bigr)=\tag 3.7\\
&\phantom{\Pi^*(\Omega')}
{}=\langle\xi_2,Iu_1\rangle
-\langle-I\xi_1,u_2\rangle=
\langle\xi_2,Iu_1\rangle
-\langle\xi_1,Iu_2\rangle,
\endalign$$
where $\xi_1,\xi_2,u_1,u_2\in{\frak m}$.
From~(3.6)
and~(3.7) it follows that the forms
$d(\Pi^*\theta')=\Pi^*(d\theta')$ and
$\Pi^*\Omega'$ coincide when restricted to
${\Cal T}_h\oplus{\Cal T}_v$ and,
consequently, on the whole tangent bundle
$T(G\times W)={\Cal K}\oplus{\Cal T}_h\oplus{\Cal T}_v$ because
${\Cal K}$ is the kernel of
$\Pi_*$. Therefore
$\Omega'=d\theta'$. Thus for the pair
$\bigl(J^-J(P),(J(P),\Omega)\bigr)$ conditions
(1) and (2) of Lemma~2.7 hold,
and we are done.
{}\hfill$\square$
\enddemo

As an immediate consequence of the proof we obtain
\proclaim{Corollary 3.3.1}
Let $P,{\bold g}$ be as in Theorem~{\rm 3.3} and
$IP=\overline PI$. Let
$\theta'$ and
$\Omega'$ be the forms given by~{\rm(3.5)}. Then
$\Omega'$ is the K\"ahler form of the K\"ahler structure
$(J^-J(P),{\bold g})$ and
$\Omega'=d\theta'$.
\endproclaim

It is evident that integrability condition~(2.7) for
$P$ is equivalent to a pair of real equations
for its real and imaginary parts $R$ and
$S$.
The following proposition establishes
more restrictive conditions for $R$ and
$S$ if the pair $\bigl(J^-,(J(P),\Omega)\bigr)$ defines a
hyper-K\"ahlerian structure.
\proclaim{Proposition 3.4}
Let
$(J(P),\Omega)$ be a K\"ahler structure on
$D$ such that $IP=\overline PI$ on $W$.
Then
$$\align
(L_\xi(R+SR^{-1}S))_w(I\eta)= \bigl[[w,IR^{-1}_w\xi],\eta\bigr],
&\quad w\in W,\ \xi,\eta\in{\frak m}; \tag 3.8\\
L_\xi(SR^{-1})\eta=L_\eta(SR^{-1})\xi,
&\quad w\in W,\ \xi,\eta\in{\frak m}. \tag 3.9\\
\endalign$$
Locally the mapping
$w\mapsto (SR^{-1})_w$ is a tangent one of some vector-function
with values in ${\frak m}$.
\endproclaim
\demo{Proof}
Put $J=J(P)$. Let
${\bold g}$ be the Hermitian metric
corresponding to $(J,\Omega)$.
By Theorem~3.3 the triple
$({\bold g},J^-,J)$ is a hyper-K\"ahlerian
structure, in particular, the 2-form
$\omega_1=\Omega^-$,
$\Omega^-(\cdot,\cdot)\overset\text{def}\to{=}
\Omega(J^-J\cdot,\cdot)$ on
$D$ is closed and so is the form
$\Pi^*\Omega^-$ on
$G\times W$. By~(2.2) the form
$(\Pi^*\Omega^-)(\cdot,\cdot)=
(\Pi^*\Omega)(\tilde J^-\tilde J\cdot,\cdot)$.
Therefore for any vector fields
$\tilde X,\tilde Y,\tilde Z$ on
$G\times W$ we have
$$
\sum_{\tilde X\tilde Y\tilde Z}
\tilde X\big((\Pi^*\Omega)(\tilde J^
-\tilde J\tilde Y,\tilde Z)\bigr)+
\sum_{\tilde X\tilde Y\tilde Z}
(\Pi^*\Omega)(\tilde J^-\tilde J\tilde X,[\tilde Y,\tilde Z])
=0
\tag 3.10$$
(we sum here over the cyclic permutations of
$\tilde X,\tilde Y,\tilde Z$).

Putting in~(3.10) $\tilde X=(0,\xi)$,
$\tilde Y=(\eta^l,0)$ and
$\tilde Z=(\chi^l,0)$ for (fixed)
$\xi,\eta,\chi\in{\frak m}$, we obtain only two non-zero terms in
the left-hand side of~(3.10) (here all objects are left
$G$-invariant) and, consequently,
$$
(L_\xi\langle -I(R+SR^{-1}S)\eta,\chi\rangle)_w- \langle
w,[-IR^{-1}\xi,[\eta,\chi]]\rangle=0,
$$
because
$\tilde J^-\tilde J(0,\xi)=((-IR^{-1}\xi)^l,-ISR^{-1}\xi)$.
Since the form $\langle ,\rangle$ is
$\operatorname{Ad} G$-invariant and $IR=RI$,
$IS=-IS$, we derive condition (3.8).

To prove (3.9), put in~(3.10)
$\tilde X=(0,\xi)$,
$\tilde Y=(0,\eta)$ and
$\tilde Z=(\chi^l,0)$.
Then we have only two non-zero terms on
the left in~(3.10), i.e. the following equality
$$
L_\xi\langle -ISR^{-1}\eta,\chi\rangle
-L_\eta\langle \xi,-IR^{-1}S\chi\rangle=0.
$$
By Proposition~2.8 the endomorphisms
$R_w$ and
$S_w$ are symmetric and by definition
$I$ is skew-symmetric (with respect to the form
$\langle ,\rangle|{\frak m}$). Using~(3.4) we obtain that
$IR^{-1}S=-R^{-1}SI$ and, consequently,
$L_\xi(SR^{-1})\eta=L_\eta(SR^{-1})\xi$. This identity then gives
the latter assertion of the proposition.
{}\hfill$\square$
\enddemo

\remark{Remark 3.5}
It is easy to verify that conditions~(3.8)
and~(3.9) are equivalent to condition (1.8)
in~\cite{BG1}.
\endremark

\proclaim{Lemma 3.6}
Let $P\in\text{Alm}(W)$ and $IP=\overline PI$. Assume that for
$P$ conditions~{\rm(3.8)} and~{\rm(3.9)}
hold. Then the almost complex structure
$J(P)$ on
$D$ is integrable iff for all
$w\in W$
$$
(SR^{-1})_w\bigl([w,[\xi,\eta]]\bigr)=
\bigl[[w,(R^{-1}S)_w\xi],\eta\bigr]
-\bigl[[w,(R^{-1}S)_w\eta],\xi\bigr],
\quad \xi,\eta\in{\frak m}.
\tag 3.11$$
\endproclaim
\demo{Proof}
Since
$R+SR^{-1}S=(1-iSR^{-1})(R+iS)$, we have
$$
L_X(R+SR^{-1}S)=-iL_X(SR^{-1})\cdot P+
(1-iSR^{-1})\cdot L_X P
$$
for
$X\in{\Gamma}(TW)$ and, consequently, by~(3.8)
and~(3.9)
$$\align
&(1-iSR^{-1})_w\Big(L_{P\xi} P(\eta)
-L_{P\eta} P(\xi)\Big)_w \\
&=i\Big((L_{P\xi}(SR^{-1}))(P\eta)
-(L_{P\eta}(SR^{-1}))(P\xi)\Big)_w\\
&+\Big(L_{P\xi}(R+SR^{-1}S)(\eta)
-L_{P\eta}(R+SR^{-1}S)(\xi)\Big)_w\\
&=\bigl[[w,I(R^{-1}P)_w\xi],-I\eta\bigr]
-\bigl[[w,I(R^{-1}P)_w\eta],-I\xi\bigr].
\endalign$$
Since the operator
$1-iSR^{-1}=\overline PR^{-1}$ is invertible and
$IP=\overline PI$, integrability condition~(2.7) holds iff
$$
(\overline PR^{-1})_w([w,[\xi,\eta]])=
\bigl[[w,(R^{-1}\overline P)_wI\xi],I\eta\bigr]
-\bigl[[w,(R^{-1}\overline P)_wI\eta],I\xi\bigr].
\tag 3.12$$
Using the Jacobi identity and properties~(3.1) of
$I$, we obtain that real part of the right-hand side of~(3.12) is
equal to $[w,[I\xi,I\eta]]=[w,[\xi,\eta]]$.
The imaginary part of~(3.12) is equivalent to~(3.11).
{}\hfill$\square$
\enddemo

\head
{\bf\S 4. Hyperk\"ahler structures on irreducible
Hermitian symmetric spaces}
\endhead

In this section the main theorem describing
K\"ahler structures anticommuting with
$J^-$ on homogeneous domains in
$T(G/K)$ is proved. We give a general formula for the K\"ahler
potential. We seek this formula in the form as in~\cite{BG1}.

We continue with the previous notations but in this section it is
assumed in addition that the Hermitian symmetric space
$G/K$ is irreducible, has compact type and the form
$\langle ,\rangle$ on
${\frak g}$ is positive-definite. In particular, then
${\frak g}$ is a simple compact Lie algebra of rank
$l$ and the subgroup $K$ is connected.

{\bf 4.1. Root theory of Hermitian symmetric spaces.}
Here we will review few facts about Hermitian symmetric
spaces~\cite{He, Ch.VIII, \S\S4--7}. The compact Lie subalgebra
${\frak k}\subset{\frak g}$ has the one-dimensional center
${\frak z}$ and coincides with the centralizer
${\frak g}^{\frak z}$ of ${\frak z}$ in
${\frak g}$, in particular,
$\text{rk\,}{\frak g}=\text{rk\,}{\frak k}=l$. Let
${\frak t}$ be some Cartan subalgebra of
${\frak k}$. Then
${\frak z}\subset{\frak t}$ and
${\frak t}$ is a Cartan subalgebra of
${\frak g}$. The complex space
${\frak t}^{\Bbb C}$ is a Cartan
subalgebra of the simple complex Lie algebra
${\frak g}^{\Bbb C}$. Let
$\Delta$ be the root system of
${\frak g}^{\Bbb C}$ with respect to
${\frak t}^{\Bbb C}$. Denote by
$\Delta_{\frak k}$ the set of roots in
$\Delta$ which vanish identically on
${\frak z}$. This is the root system of
$({\frak k}^{\Bbb C},{\frak t}^{\Bbb C})$. Put
$\Delta_{\frak m}=\Delta\setminus\Delta_{\frak k}$. Then we have
the direct decompositions
$$
{\frak g}^{\Bbb C}={\frak t}^{\Bbb C}\oplus\sum_{\alpha\in\Delta}
\tilde{\frak g}^\alpha,
\qquad
{\frak k}^{\Bbb C}={\frak t}^{\Bbb C}\oplus
\sum_{\alpha\in\Delta_{\frak k}}\tilde{\frak g}^\alpha,
\qquad
{\frak m}^{\Bbb C}=\sum_{\alpha\in\Delta_{\frak m}}\tilde{\frak g}^\alpha.
$$
The algebra
${\frak k}$ is a maximal subalgebra of
${\frak g}$~\cite{GG, (8.5.1)}. Since this subalgebra is
$\operatorname{ad}{\frak t}$-invariant, according to~\cite{GG,
(8.3.7)} there exists a system of simple roots
$\pi=\{\alpha_1,..,\alpha_l\}\subset\Delta$ such that for any
$\alpha\in\Delta$ we have
$n_1^\alpha\in\{0,1,-1\}$, where
$\alpha=n_1^\alpha\alpha_1+\cdots +n_l^\alpha\alpha_l$. Then
$\Delta_{\frak k}=\{\alpha\in\Delta: n_1^\alpha=0\}$. Denote by
$\Delta^+\subset\Delta$ the corresponding
set of positive roots.

Choose for any
$\alpha\in\Delta^+$ a triple
$(H_\alpha,E_\alpha,E_{-\alpha})\in
i{\frak t}\times\tilde{\frak g}^\alpha
\times\tilde{\frak g}^{-\alpha}$
such that
$[H_\alpha,E_{\pm\alpha}]=\pm 2E_{\pm\alpha}$ and
$[E_{\alpha},E_{-\alpha}]=-H_\alpha$.
This choice can be made so that
${\frak g}$ has a basis consisting of a basis of
${\frak t}$ and
$X_\alpha=\frac12(E_\alpha+E_{-\alpha})$,
$Y_\alpha=\frac i2(E_\alpha-E_{-\alpha})$,
$\alpha\in\Delta^+$ (the space
${\frak t}$ is spanned by the vectors
$T_\alpha=\frac i2H_\alpha$,
$\alpha\in\Delta^+$). One has
$$
[X_\alpha,Y_\alpha]=T_\alpha, \quad
[T_\alpha,X_\alpha]=Y_\alpha, \quad
[T_\alpha,Y_\alpha]=-X_\alpha,
\quad\alpha\in\Delta^+.
$$
Putting
$\Delta^+_{\frak k}=\Delta_{\frak k}\cap\Delta^+$ and
$\Delta^+_{\frak m}=\Delta_{\frak m}\cap\Delta^+$, we obtain
$$
{\frak k}={\frak t}\oplus\sum_{\alpha\in\Delta^+_{\frak k}}
({\Bbb R} X_\alpha\oplus{\Bbb R} Y_\alpha),
\qquad
{\frak m}=\sum_{\alpha\in\Delta^+_{\frak m}}
({\Bbb R} X_\alpha\oplus{\Bbb R} Y_\alpha).
$$
Since for each pair
$\alpha,\beta\in\Delta^+_{\frak m}$
the sum $\alpha+\beta$ is not a root,
$-2i\alpha(T_\beta)=\alpha(H_\beta)\ge0$.
There exists a unique element
$Z_0\in{\frak z}$ such that
$\alpha_1(Z_0)=i$
($\alpha_j|{\frak z}=0$ for $j=\overline{2,l}$). Putting
$I=\operatorname{ad}_{Z_0}|{\frak m}:{\frak m}\to{\frak m}$,
we see that $IX_\alpha=Y_\alpha$ and
$IY_\alpha=-X_\alpha$ for all
$\alpha\in\Delta^+_{\frak m}$. Then
$$
[T,\xi_\alpha]=-i\alpha(T)I\xi_\alpha,
\qquad\text{where}\
T\in{\frak t},\ \xi_\alpha\in
({\Bbb R} X_\alpha\oplus{\Bbb R} Y_\alpha)\subset{\frak m}.
\tag 4.1$$
Hence
$I^2=-\operatorname{Id}_{\frak m}$. Moreover,
$I$ and the automorphism
$\exp(\frac\pi2\operatorname{ad}_{Z_0})\in\operatorname{Ad}(K)$
coincide when restricted to
${\frak m}$. Since the Lie group
$K$ is connected, the group
$\operatorname{Ad}(K)$ commutes elementwise with
$\exp(\frac\pi2\operatorname{ad}_{Z_0})$.

Whereas all Cartan subspaces (maximal abelian subalgebras) of
${\frak m}$ are conjugate under the linear isotropy group
$\operatorname{Ad}(K)$ it is possible in the special
situation here to select such a Cartan subspace
${\frak a}$ with particular reference to
$\Delta$. Two roots
$\alpha,\beta\in\Delta$ are called strongly orthogonal if
$\alpha\pm\beta\not\in(\Delta\cup\{0\})$. There exists a subset of
$\Delta^+_{\frak m}$ consisting of
$r=\text{rk\,}(G/K)$ strongly orthogonal roots
$\beta_1,..,\beta_r$~\cite{He, Ch.VIII,
Prop.7.4}. Then the subspaces
${\frak a}=\sum_{j=1}^r{\Bbb R} X_{\beta_j}$ and
$I{\frak a}=\sum_{j=1}^r{\Bbb R} Y_{\beta_j}$
are Cartan subspaces of
${\frak m}$; the Lie subalgebra of
${\frak g}$ generated by subspaces
${\frak a}$ and
$I{\frak a}$ is isomorphic to the semisimple compact Lie algebra
${\hat{\frak g}}=\bigoplus_{j=1}^r{\hat{\frak g}}_j$, where each
${\hat{\frak g}}_j=({\Bbb R} X_{\beta_j}\oplus
{\Bbb R} Y_{\beta_j} \oplus {\Bbb R} T_{\beta_j})\simeq su(2)$.
Then
$-i\beta_k(T_{\beta_j})=\frac12\beta_k(H_{\beta_j})=\delta^k_j$.
We have
\proclaim{Proposition 4.1}~{\rm\cite{He}}
Any Cartan subspace ${\frak a}$ of
${\frak m}$ has the form
${\frak a}=\sum_{j=1}^{r}{\Bbb R} X_{\beta_j}$.
The Lie subalgebra of
${\frak g}$ generated by subspaces ${\frak a}$ and
$I{\frak a}$ is isomorphic to the semisimple compact Lie algebra
${\hat{\frak g}}=\bigoplus_{j=1}^r{\hat{\frak g}}_j$, where each
${\hat{\frak g}}_j\simeq su(2)$.
\endproclaim

Denote by
${\frak g}^{\frak a}$ and
${\frak k}^{\frak a}$ the centralizers of the Cartan subspace
${\frak a}$ in ${\frak g}$ and
${\frak k}$ respectively. By~(2.3),
${\frak g}^{\frak a}={\frak a}\oplus{\frak k}^{\frak a}$. In particular,
$\text{rk\,}{\frak k}^{\frak a}=\text{rk\,}{\frak g}-r$ and
by~(3.1)
$[I{\frak a},{\frak k}^{\frak a}]=I[{\frak a},{\frak k}^{\frak a}]=0$.
\proclaim{Proposition 4.2}
Let
${\frak a}\subset{\frak m}$ be a Cartan
subspace. For the irreducible symmetric space
$G/K$ either
${\frak a}$ is a Cartan subalgebra of
${\frak g}$ {\rm(}and
${\frak k}^{\frak a}=0${\rm)} or the centralizer of the algebra
${\frak k}^{\frak a}$ in
${\frak m}$ coincides with the space
${\frak a}\oplus I{\frak a}$.
\endproclaim
\demo{Proof}
Let
${\frak n}({\frak k}^{\frak a})$ be the normalizer of the algebra
${\frak k}^{\frak a}$ in
${\frak g}$. Since
${\frak k}^{\frak a}$ is an ideal of the compact algebra
${\frak n}({\frak k}^{\frak a})$,
we obtain the following splitting
${\frak n}({\frak k}^{\frak a})=
{\frak k}^{\frak a}\oplus{\frak g}_*$,
where
${\frak g}_*=\{X\in{\frak g}: [X,{\frak k}^{\frak a}]=0,
\langle X,\ {\frak k}^{\frak a}\rangle=0\}$.
It is clear that the rank of the Lie algebra
${\frak g}_*$ does not exceed
$(\text{rk\,}{\frak g}-\text{rk\,}{\frak k}^{\frak a})=r$.
By the Jacobi identity
$[[{\frak a},I{\frak a}],{\frak k}^{\frak a}]=0$. Therefore
$[{\hat{\frak g}},{\frak k}^{\frak a}]=0$ and
$\langle {\hat{\frak g}},{\frak k}^{\frak a}\rangle=0$. Hence
${\frak g}_*$ is a semisimple algebra of rank
$r$ containing the semisimple subalgebra
${\hat{\frak g}}$ of maximal rank. Since
$Z_0$ is an element of the center of the subalgebra
${\frak k}$, we obtain that
$\operatorname{ad}_{Z_0}({\frak g}_*)\subset{\frak g}_*$
and the triple
$({\frak g}_*,{\frak k}_*,I_*)$, where
${\frak k}_*={\frak k}\cap{\frak g}_*$,
${\frak m}_*={\frak m}\cap{\frak g}_*$,
$I_*=I|{\frak m}_*$, is an Hermitian orthogonal
symmetric Lie algebra. Taking into account that
${\frak a}\subset{\frak m}_*$, we conclude that
this orthogonal Lie algebra has maximal possible rank
$r=\text{rk\,}{\frak g}_*$. Its each irreducible component also
has maximal possible rank
(${\frak a}$ is a Cartan subalgebra of ${\frak g}_*$),
i.e. corresponds to the compact symmetric space
$Sp(n)/U(n)$ for appropriate
$n\ge1$~\cite{He, Ch.X, \S 6}. Therefore
${\frak g}_*$ is a semisimple Lie algebra of type
$C_{n_1}\oplus C_{n_2}\oplus\cdots$. We claim that
$n_1=n_2=..=n_r=1$ or
${\frak k}^{\frak a}=0$.

Indeed, by construction the algebra
${\frak g}_*\oplus{\frak k}^{\frak a}$ is a subalgebra of
${\frak g}$ of maximal rank. Hence
${\frak g}_*$ is a regular subalgebra of
${\frak g}$, i.e.
$[{\frak t}_1,{\frak g}_*]\subset{\frak g}_*$
for some Cartan subalgebra
${\frak t}_1$ of
${\frak g}$. But the algebra
${\frak g}$ is a compact Lie algebra from the following list
$A_l,B_l,C_l,D_l,E_6,E_7$~\cite{He,
Ch.X, \S 6}. Since for algebras
$A_l,D_l,E_6,E_7$ all roots of their root systems have the same
length, these algebras do not contain regular subalgebras of type
$C_n,n\ge2$. So we have to consider only two cases when
$G/K$ is the symmetric space
$SO(2l+1)/(SO(2)\times SO(2l-1))$ with $l\ge3$ or
$Sp(l)/U(l)$ with $l\ge2$. In the first case
${\frak k}^{\frak a}\simeq so(2l-3)$ and
${\frak g}_*\simeq so(4)\simeq C_1\oplus C_1$~\cite{He,
Ch.X}. In the second case ranks of $G$ and
$G/K$ coincide, i.e.
${\frak k}^{\frak a}=0$. So the claim is proved.
{}\hfill$\square$
\enddemo

{\bf 4.2. Invariant mappings and root
theory of Hermitian symmetric spaces.}
In this subsection the transformation
of the restricted root system of
$({\frak g},{\frak a})$ induced by the action of
$I$ on ${\frak m}$ is studied.

Let ${\frak a}$ be some Cartan subspace of
${\frak m}$.
For each
$\lambda$ in the dual space of
${\frak a}^{\Bbb C}$ let
$\tilde{\frak g}_\lambda=\{X\in{\frak g}^{\Bbb C}:
[H,X]=\lambda(H), H\in{\frak a}^{\Bbb C}\}$.
Then
$\lambda$ is called a restricted root if
$\lambda\ne0$ and
$\tilde{\frak g}_\lambda\ne0$.
The set of all such $\lambda$ is denoted by $\Sigma$.
The simultaneous diagonalization of
$\operatorname{ad}({\frak a}^{\Bbb C})$ in
${\frak g}^{\Bbb C}$ gives the decomposition
${\frak g}^{\Bbb C}=\tilde{\frak g}_0\oplus\sum_{\lambda\in\Sigma^+}
(\tilde{\frak g}_\lambda\oplus\tilde{\frak g}_{-\lambda})$,
where
$\Sigma^+$ is an arbitrary subset of positive restricted roots
in $\Sigma$.

For each linear form $\lambda$ on ${\frak a}^{\Bbb C}$ put
$$\align
{\frak m}_\lambda&\overset\text{def}\to{=} \{\eta\in{\frak m}:
\operatorname{ad}^2_w(\eta)
=\lambda^2(w)\eta,\ \forall w\in{\frak a}\}, \\
{\frak k}_\lambda&\overset\text{def}\to{=} \{\zeta\in{\frak k}:
\operatorname{ad}^2_w(\zeta)=\lambda^2(w)\zeta,\
\forall w\in{\frak a}\}.
\endalign$$
Then ${\frak m}_{\lambda}={\frak m}_{-\lambda}$,
${\frak k}_{\lambda}={\frak k}_{-\lambda}$,
${\frak m}_0={\frak a}$ and ${\frak k}_0$ equals
${\frak k}^{\frak a}$, the centralizer of ${\frak a}$ in
${\frak k}$. By~\cite{He, Lemma 11.3, Ch.VII} the following
decompositions are direct and orthogonal:
$$
{\frak m}={\frak a}\oplus\sum_{\lambda\in\Sigma^+}{\frak m}_\lambda,
\qquad
{\frak k}={\frak k}^{\frak a}\oplus
\sum_{\lambda\in\Sigma^+}{\frak k}_\lambda.
\tag 4.2$$
We need the following lemma which is a weak
generalization of Lemma 2.3 in~\cite{He, Ch.VII}.
\proclaim{Lemma 4.3}
For any vector
$\xi_\lambda\in{\frak m}_\lambda$,
$\lambda\in\Sigma^+$ there exists a unique vector
$\zeta_\lambda\in{\frak k}_\lambda$ such that
$$
[w,\xi_\lambda]=i\lambda(w)\zeta_\lambda,
\quad [w,\zeta_\lambda]=-i\lambda(w)\xi_\lambda
\qquad\text{for all}\ w\in{\frak a}.
\tag 4.3$$
In particular, $\dim{\frak m}_\lambda=\dim{\frak k}_\lambda$ and
$$
\operatorname{ad}_{w'}\operatorname{ad}_{w''}(\xi_\lambda)
=\lambda(w')\lambda(w'')\xi_\lambda,
\quad\text{where}\quad w',w''\in{\frak a}.
$$
\endproclaim
\demo{Proof}
For completeness and mainly to fix the notation
we shall prove this lemma. It is clear that
$({\frak m}_\lambda\oplus{\frak k}_\lambda)^{\Bbb C}
=(\tilde{\frak g}_\lambda\oplus\tilde{\frak g}_{-\lambda})$.
Therefore
$\operatorname{ad}_w({\frak m}_\lambda)\subset{\frak k}_\lambda$
and
$\operatorname{ad}_w({\frak k}_\lambda)\subset{\frak m}_\lambda$
for $w\in{\frak a}$. So the endomorphisms
$\operatorname{ad}_w$ and
$\operatorname{ad}^2_w$ when restricted to
${\frak m}_\lambda\oplus{\frak k}_\lambda$ are nondegenerate or
degenerate simultaneously. Hence the subspace
$\{[w,\xi_\lambda],\ w\in{\frak a}\}\subset{\frak k}_\lambda$
is one-dimensional. Since
$\lambda({\frak a})\in i{\Bbb R}$, there is the element
$\zeta_\lambda\in{\frak k}_\lambda$ such that for the pair
$\{\xi_\lambda,\zeta_\lambda\}$ condition~(4.3) holds. Now the latter
assertion of the lemma is evident.
{}\hfill$\square$
\enddemo

Let $f:W\to{\frak m}$ be a mapping.
Identifying the tangent spaces
$T_w{\frak m}$ and $T_{f(w)}{\frak m}$ with
${\frak m}$, we can consider the tangent mapping
$f_{*w}:T_w{\frak m}\to T_{f(w)}{\frak m}$ as an endomorphism on
${\frak m}$. We say
$f$ is $K$-equivariant if
$\operatorname{Ad}_k\circ f=f\circ\operatorname{Ad}_k$ on $W$ for all
$k\in K$. For such a mapping its tangent map
$f_*: w\mapsto f_{*w}\in\text{End}({\frak m})$ is also
$K$-equivariant, i.e. satisfies~(2.6).
Denote by $\text{EC}(W)$ the set of all
$K$-equivariant mappings $f$ (on $W$) which
leave some (and, consequently, each) Cartan subspace ${\frak a}$
of ${\frak m}$ invariant, i.e. $f(W\cap{\frak a})\subset{\frak a}$.

Let $f\in\text{EC}(W)$. By $K$-equivariance,
$f_{*w}([\zeta,w])=[\zeta,f(w)]$ for any
$\zeta\in{\frak k}$. Hence by~(4.3) for any
$w\in W\cap{\frak a}$ and
$\xi_\lambda\in{\frak m}_\lambda$, $\lambda\in\Sigma^+$
$$
f_{*w}(i\lambda(w)\xi_\lambda)=
f_{*w}([\zeta_\lambda,w])= [\zeta_\lambda,f(w)]=
i\lambda(f(w))\xi_\lambda.
\tag 4.4$$

For future use we next prove two lemmas.
\proclaim{Lemma 4.4}
Let $f\in\text{EC}(W)$. If
$w\in W\cap{\frak a}$ then
$f_{*w}({\frak a})\subset{\frak a}$ and for each
$\lambda\in\Sigma^+$ the subspace
${\frak m}_\lambda$ is an eigenspace of
$f_{*w}$ with the eigenvalue
$\lambda(f(w))/\lambda(w)$.
Moreover, for all $w\in W$, $\xi,\eta\in{\frak m}$
$$
f_{*w}\bigl([w,[\xi,\eta]]\bigr)=
\bigl[[w,f_{*w}(\xi)],\eta\bigr]
-\bigl[[w,f_{*w}(\eta)],\xi\bigr].
\tag 4.5$$
\endproclaim
\demo{Proof}
Let $w\in W\cap{\frak a}$.
Suppose that
$\xi\in{\frak a}$. Since by definition
$f_{*w}({\frak a})\subset{\frak a}$, the first term
on the right in~(4.5) vanishes. But
$[\operatorname{ad}_w,\operatorname{ad}_\xi]=0$.
Therefore~(4.5) is equivalent to the relation
$$
[(\operatorname{ad}_{w}
\operatorname{ad}_\xi)|{\frak m},f_{*w}]=0,
$$
which holds, because
$\operatorname{ad}_{w}\operatorname{ad}_\xi({\frak a})=0$
and each subspace ${\frak m}_\lambda$,
$\lambda\in\Sigma^+$ is an eigenspace of endomorphisms
$\operatorname{ad}_{w}\operatorname{ad}_\xi$ and
$f_{*w}$.

Since relation~(4.5) is
skew-symmetric for exchanges of two variables
$\xi$ and
$\eta$, it remains to prove~(4.5) if
$\xi=\xi_\lambda\in{\frak m}_\lambda$ and
$\eta=\xi'_\nu\in{\frak m}_\nu$ for
$\lambda,\nu\in\Sigma^+$. There are vectors
$\zeta_\lambda\in{\frak k}_\lambda$ and
$\zeta'_\nu\in{\frak k}_\nu$ such that for the pairs
$(\xi_\lambda,\zeta_\lambda)$ and
$(\xi'_\nu,\zeta'_\nu)$ condition~(4.3) holds for all
$w\in W\cap{\frak a}$. But
$$
[{\frak k}_\lambda,{\frak m}_\nu]+[{\frak m}_\lambda,{\frak k}_\nu]
\subset{\frak m}_{\lambda+\nu}+{\frak m}_{\lambda-\nu}
\qquad\text{and}\qquad
[{\frak m}_\lambda,{\frak m}_\nu]
\subset{\frak k}_{\lambda+\nu}+{\frak k}_{\lambda-\nu}
$$
(see \cite{He, Ch.VII, Lemma 11.4}). In particular,
$[\xi_\lambda,\xi'_\nu]=\zeta_+ + \zeta_-$, where
$\zeta_{\pm}\in{\frak k}_{\lambda\pm\nu}$.
If $\lambda-\nu=0$ then $[{\frak a},\zeta_-]=0$. Therefore
there exist vectors
$\xi_\pm\in{\frak m}_{\lambda\pm\nu}$ such that
$[w',\zeta_\pm]=-i(\lambda\pm\nu)(w')\xi_\pm$ for all
$w'\in{\frak a}$. Now from the Jacobi identity
$$
[w',[\xi_\lambda,\xi'_\nu]]
=[[w',\xi_\lambda],\xi'_\nu]
-[[w',\xi'_\nu],\xi_\lambda]
$$
it follows that
$$
-i(\lambda+\nu)(w')\xi_+ -i(\lambda-\nu)(w')\xi_-
=i\lambda(w')[\zeta_\lambda,\xi'_\nu]
-i\nu(w')[\zeta'_\nu,\xi_\lambda].
\tag 4.6$$
Taking into account that
$f_{*w}(\xi_\pm)={\displaystyle\frac{(\lambda\pm\nu)(f(w))}
{(\lambda\pm\nu)(w)}}\xi_\pm$
and similar relations hold for
$\xi_\lambda$,
$\xi'_\nu$, we obtain~(4.5) replacing
$w'$ by
$f(w)$ in identity~(4.6). Noting that ${\frak m}$ is a union of its
Cartan subspaces, we complete the proof.
{}\hfill$\square$
\enddemo

\proclaim{Lemma 4.5}
Let $f\in\text{EC}(W)$ and let
$\sigma$ be a 1-form on $W$ such that
$\sigma_w(\cdot)=\langle f(w),\cdot\rangle$.
Then
\roster
\item
the form $\sigma$ is $\operatorname{Ad}(K)$-invariant;
\item
$\sigma$ is closed iff so is its restriction to the set
$W\cap{\frak a}$.
\endroster
\endproclaim
\demo{Proof}
It is immediate that $\sigma$
is invariant. By definition
$$
d\sigma_w(\xi,\eta)=\langle f_{*w}(\xi),\eta \rangle
-\langle f_{*w}(\eta),\xi \rangle,
\quad \xi,\eta\in{\frak m}.
$$
For $w\in{\frak a}$, by Lemma~4.4
$f_{*w}({\frak a})\subset{\frak a}$,
$f_{*w}({\frak a}^\bot)\subset{\frak a}^\bot$ and the restriction
$f_{*w}|{\frak a}^\bot$ is a symmetric operator.
Therefore $d\sigma_w=0$ iff
$d\sigma_w({\frak a},{\frak a})=0$ for all
$w\in W\cap{\frak a}$.
{}\hfill$\square$
\enddemo

Let
$O\subset{\Bbb R}$ be a domain containing spectrums of all operators
$(-\operatorname{ad}_w^2-\operatorname{ad}_{Iw}^2)|{\frak m}$,
$w\in W$. For a real-analytic function
$q$ on $O$, define a mapping
$\hat q:W\to{\frak m}$ by
$\hat q(w)=q(-\operatorname{ad}_w^2-\operatorname{ad}_{Iw}^2)(w)$.

\proclaim{Proposition 4.6}
The mapping $\hat q:W\to{\frak m}$ is
$K$-equivariant and each Cartan subspace
${\frak a}$ of
${\frak m}$ is invariant with respect to
$\hat q$, i.e.
$\hat q\in\text{EC}(W)$. Moreover,
$$
{\frak m}=\sum_{j=1}^r {\Bbb R} X_j\oplus
\sum_{\lambda\in\Sigma^+}{\frak m}_\lambda,
\qquad\text{where}\ X_j=X_{\beta_j},
\tag 4.7$$
is the orthogonal eigenspace splitting for all
$\hat q_{*w}$, $w\in W\cap{\frak a}$.
\endproclaim
\demo{Proof}
Since the endomorphism
$I$ belongs to the center of the group $\operatorname{Ad}(K)|{\frak m}$,
it follows that
$\operatorname{Ad}_k\circ\operatorname{ad}_{Iw}
=\operatorname{ad}_{I(\operatorname{Ad}_k w)}\circ\operatorname{Ad}_k$
and, consequently,
$\operatorname{Ad}_k\circ\hat q=\hat q\circ\operatorname{Ad}_k$ on
$W$ for all
$k\in K$. Now fix some Cartan subspace
${\frak a}=\sum_{j=1}^{r}{\Bbb R} X_{\beta_j}$ of
${\frak m}$ (as in subsection 4.1) and relabel
$X_{\beta_j},Y_{\beta_j},T_{\beta_j}$
to read $X_j,Y_j,T_j$, $j=\overline{1,r}$.
One has
$$
[X_j,Y_k]=\delta_j^k T_j, \quad
[T_j,X_k]=\delta_j^k Y_j, \quad
[T_j,Y_k]=-\delta_j^k X_j,\quad j,k=\overline{1,r}.
\tag 4.8$$
In particular, $IX_j=Y_j$ and
$IY_j=-X_j$. Since
$\operatorname{ad}^2_{IX_j}(X_k)=-\delta^j_k X_k$
and $[{\frak a},{\frak a}]=0$,
for any $w=\sum_{j=1}^{r}x_j X_j\in W\cap{\frak a}$ we have
$$
\hat q(w)=\sum_{j=1}^{r}x_j q(x_j^2)X_j
\quad\text{and}\quad
\hat q_{*w}(X_j)=\bigl(q(x_j^2)+2x_j^2q'(x_j^2)\bigr)X_j,\
j=\overline{1,r},
\tag 4.9$$
i.e. $\hat q(W\cap{\frak a})\subset{\frak a}$.
The latter assertion follows immediately
from~(4.4) and~(4.9).
{}\hfill$\square$
\enddemo

Fix in the Cartan subspace
${\frak a}\subset{\frak m}$ a basis
$\{X_j\}_{j=1}^r$~(4.8). Let the set of restricted roots
$\Sigma$ of
$({\frak g},{\frak a})$ be ordered
lexicographically with respect to the basis
$\{-iX_j\}_{j=1}^r$ in
$i{\frak a}\subset{\frak a}^{\Bbb C}$ (all
$\lambda\in\Sigma$ are real on the subspace
$i{\frak a}\subset{\frak g}^{\Bbb C}$). Denote by
$\Sigma^+$ the corresponding system of
positive restricted roots. Choose the basis
$\{\epsilon_j\}_{j=1}^{r}$ in the (complex) space
$({\frak a}^{\Bbb C})^*$ dual to the basis
$\{X_j\}_{j=1}^{r}$ of ${\frak a}^{\Bbb C}$.
\proclaim{Proposition 4.7}
For any vector
$\xi\in{\frak m}_\lambda$,
$\lambda\in\Sigma^+$ the vector
$I\xi$ belongs to the subspace
${\frak m}_{\lambda_I}$,
$\lambda_I\in\Sigma^+\cup\{0\}$, i.e.
$\operatorname{ad}_w^2(\xi)=\lambda^2(w)\xi$,
$\forall w\in{\frak a}$ implies
$\operatorname{ad}_w^2(I\xi)=\lambda_I^2(w)(I\xi)$.
The set $\Sigma^+$ is a subset of the set
$$
(BC)_r^+=
\left\{
{\tfrac i2}\epsilon_j,\ i\epsilon_j,\
j=\overline{1,r};\quad
{\tfrac i2}(\epsilon_p\pm\epsilon_k),\
1\le p< k\le r\right\},
$$
and the set
$\{(\lambda,\lambda_I),\lambda\in\Sigma^+\}$
is a subset of the set
$$
\left\{
\bigl({\tfrac i2}\epsilon_j,
{\tfrac i2}\epsilon_j\bigr),
\bigl(i\epsilon_j,0\bigr),
\ j=\overline{1,r};
\bigl({\tfrac i2}
(\epsilon_p\pm\epsilon_k),
{\tfrac i2}
(\epsilon_p\mp\epsilon_k)\bigr),\
1\le p<k\le r
\right\}.
$$
\endproclaim
\demo{Proof}
Fix $\lambda\in\Sigma^+$. Then for any vector
$w=\sum_{j=1}^r x_jX_j\in{\frak a}$:
$\lambda(w)=\sum_{j=1}^r ic_jx_j$, where
$c_j\in{\Bbb R}$. Applying~(4.4) and~(4.9)  to the function
$\hat q$ with $q(z)=z^{n}$,
$n\in{\Bbb N}$, we obtain that if
$\lambda(x_1,..,x_r)=0$ then
$\lambda(x_1^{2n+1},..,x_r^{2n+1})=0$. Therefore
$$
\left(
\matrix
x_1&x_2&..&x_r \\ x_1^3&x_2^3&..&x_r^3 \\
..&..&..&..\\
x_1^{2r-1}&x_2^{2r-1}&..&x_r^{2r-1} \\
\endmatrix
\right) \left(
\matrix
c_1\\
c_2\\
.. \\
c_r \\
\endmatrix
\right)
=0
$$
for all real vectors $(x_1,..,x_r)\in\ker\lambda$.
Since the column $(c_1,..,c_r)$ is nonzero, the determinant
$\pm\prod\limits_{j=1}^r x_j
\prod\limits_{1\le p<j\le r} (x_p^2-x_k^2)$
of the matrix
above $=0$ at these points. In other words,
$\ker\lambda=\bigcup_{\sigma\in(BC)_r^+}(\ker\lambda\cap\ker\sigma)$.
Hence
$\ker\lambda=\ker\sigma$ for some $\sigma\in(BC)_r^+$
because $(BC)_r^+$ is a finite set,  or
equivalently, $a_\lambda^{-1}\lambda\in(BC)_r^+$
for some $a_\lambda>0$.

Now put $q(z)=-z$. Then
$\hat q(w)=[Iw,[Iw,w]]$ and by~(3.1)
$$\align
\hat q_{*w}(\eta)&\overset\text{def}\to{=}
[I\eta,[Iw,w]]+[Iw,[I\eta,w]]+[Iw,[Iw,\eta]] \\
&=[I\eta,[Iw,w]]-2I[w,[w,I\eta]]
=\operatorname{ad}_{[w,Iw]}(I\eta)-2(I\operatorname{ad}_w^2)(I\eta)
\endalign$$
for each $\eta\in{\frak m}$.
Using the Jacobi identity for the vectors $I\eta,Iw,w$ and
relations~(3.1) again, we calculate the vector
$\operatorname{ad}_{[w,Iw]}(I\eta)$:
$$
[[w,Iw],I\eta]=
[w,[Iw,I\eta]]-[Iw,[w,I\eta]]=
[w,[w,\eta]]-I[w,[w,I\eta]].
\tag 4.10$$
So that
$\hat q_{*w}(\eta)=\operatorname{ad}^2_{w}(\eta)
-3(I\operatorname{ad}^2_w)(I\eta)$. Now
from~(4.4) for given
$\xi_\lambda\in{\frak m}_\lambda$
it follows that
$$
\lambda^2(x_1,..,x_r)\xi_\lambda
-3(I\operatorname{ad}^2_w)(I\xi_\lambda)=
\frac{\lambda(-x_1^3,..,-x_r^3)}
{\lambda(x_1,..,x_r)}\xi_\lambda.
\tag 4.11$$
Applying
$I$ to equation~(4.11) we obtain that
$I\xi_\lambda$ is a common eigenvector of all endomorphisms
$\operatorname{ad}^2_w$,
$w\in{\frak a}$. Hence there is a unique element
$\lambda_I\in\Sigma^+\cup\{0\}$ such that
$\operatorname{ad}^2_{w}(I\xi_\lambda)=\lambda_I^2(w)I\xi_\lambda$ and
$$
\lambda^2(x_1,..,x_r)+3\lambda_I^2(x_1,..,x_r)=
-\frac{\lambda(x_1^3,..,x_r^3)}{\lambda(x_1,..,x_r)}.
\tag 4.12$$
It is easy to verify that
if $\lambda=a_\lambda\cdot\tfrac i2(\epsilon_p\pm\epsilon_k)$, $p\ne k$,
then the pair
$(\lambda,\lambda_I)$ satisfies equation~(4.12) iff
$a_\lambda=1$ and $\lambda_I=\tfrac i2(\epsilon_p\mp\epsilon_k)$.

Since
$\operatorname{ad}^2_{X_j}(IX_k)=-\delta^j_k(IX_k)$, the covectors
$i\epsilon_j$,
$j=\overline{1,r}$ are positive restricted roots from
$\Sigma^+$. Therefore if the restricted roots
$\sigma,\lambda\in\Sigma^+$ are proportional, then
$\sigma=i\epsilon_j$ for some
$j\in\{1,..,r\}$ and $\lambda$ equals
$\tfrac i2\epsilon_j$ or
$2i\epsilon_j$. In this case all possible solutions
$(\lambda,\lambda_I)$,
$\lambda\in\Sigma^+$ of~(4.12) are pairs
$(i\epsilon_j,0)$ and
$(\tfrac i2\epsilon_j,\tfrac i2\epsilon_j)$,
$j=\overline{1,r}$.
{}\hfill$\square$
\enddemo

Let
$$\Sigma^{++}=\Sigma^+\cap
\left\{
{\tfrac i2}\epsilon_j,\ i\epsilon_j,\
j=\overline{1,r};\
{\tfrac i2}(\epsilon_p+\epsilon_k),\
1\le p<k\le r
\right\}.
$$
So
$$\align
{\frak m}=\sum_{\lambda\in\Sigma^{++}}{\frak M}_\lambda,
\qquad\text{where}\quad
{\frak M}_\lambda={\frak m}_\lambda+I{\frak m}_\lambda,
\tag 4.13\\
\endalign$$
is an orthogonal splitting of ${\frak m}$.
Since ${\frak a}={\frak m}_0$ and $\dim{\frak a}=r$, we have
\proclaim{Corollary 4.7.1}
The covector $\lambda=i\epsilon_j$,
$1\le j\le r$ is a positive restricted root from
$\Sigma^+$ with multiplicity one and
${\frak m}_\lambda={\Bbb R}(IX_j)$. If $\lambda\in\Sigma^+
\setminus\{i\epsilon_1,..,i\epsilon_r\}$ then
$I{\frak m}_\lambda={\frak m}_{\lambda_I}$.
\endproclaim

Let ${\frak t}'$ be the subspace of
${\frak t}$ spanned by
$\{T_j\}_{j=1}^{r}$. Choose the basis
$\{\epsilon'_j\}_{j=1}^{r}$ in the (complex) space
$({{\frak t}'}^{\Bbb C})^*$ dual to the basis
$\{T_j\}_{j=1}^{r}$ of ${{\frak t}'}^{\Bbb C}$. Let
$\rho$ denote the restriction mapping
$({\frak t}^{\Bbb C})^*\to({{\frak t}'}^{\Bbb C})^*$ and let
$\rho_{\frak m}$ be the mapping from
$(BC)_r^+\subset({\frak a}^{\Bbb C})^*$
to $({{\frak t}'}^{\Bbb C})^*$ with the graph given by
$$
\bigl\{
({\tfrac i2}\epsilon_j,\
{\tfrac i2}\epsilon'_j),
(i\epsilon_j,i\epsilon'_j),\ j=\overline{1,r};\
\bigl({\tfrac i2}(\epsilon_p\pm\epsilon_k),
{\tfrac i2}(\epsilon'_p+\epsilon'_k)\bigr),\ 1\le p<k\le r
\bigr\}.
\tag 4.14$$

\proclaim{Corollary 4.7.2}
For each $j=\overline{1,r}$, $\rho(\beta_j)=i\epsilon'_j$.
The set $\rho(\Delta^+_{\frak m})$ is a subset of the set
$$
\left\{
{\tfrac i2}\epsilon'_j,\ i\epsilon'_j,\
j=\overline{1,r};\quad
{\tfrac i2}(\epsilon'_p+\epsilon'_k),\
1\le p<k\le r
\right\}
\tag 4.15$$
For each $\lambda\in\Sigma^+$,
$$
{\frak m}_\lambda + I{\frak m}_\lambda=
\sum_{\alpha\in\rho^{-1}(\rho_{\frak m}(\lambda))}
({\Bbb R} X_\alpha\oplus{\Bbb R} Y_\alpha).
$$
In particular, the space ${\frak m}_\lambda+I{\frak m}_\lambda$
is an eigenspace of $I\operatorname{ad}_T$, $T\in{\frak t}'$
with the eigenvalue $i\rho_{\frak m}(\lambda)(T)$.
\endproclaim
\demo{Proof}
Choose $\lambda\in\Sigma^{+}$ and
$\xi_\lambda\in{\frak m}_\lambda$. By~(4.8)
$T_j=[X_j,IX_j]$. Putting in relation~(4.10)
$w=X_j$ and
$\eta=-I\xi_\lambda$ we obtain that
$\operatorname{ad}_{T_j}(\xi_\lambda)
=-(\lambda_I^2(X_j)+\lambda^2(X_j))I\xi_\lambda$.
But $[I,\operatorname{ad}_{T_j}|{\frak m}]=0$. Therefore
$$
[T,\xi]=-\Bigl(\sum_{j=1}^{r}
\bigl(\lambda^2(X_j)+\lambda_I^2(X_j)\bigr)
\epsilon'_j(T)\Bigr)\cdot I\xi,
\quad
\forall T\in{\frak t}',\
\xi\in({\frak m}_\lambda+I{\frak m}_{\lambda}).
\tag 4.16$$
Thus $\xi\pm iI\xi\in{\frak m}^{\Bbb C}$ are the root vectors
corresponding to the roots
$\pm i\sum_{j=1}^{r}
(\lambda^2(X_j)+\lambda_I^2(X_j))\epsilon'_j$
from
$\rho(\Delta_{\frak m})$. Now from splitting~(4.13) and
Proposition~4.7 it follows that
$\rho(\Delta^+_{\frak m})\subset Q\cup(-Q)$, where
$Q$ is set~(4.15).
Taking into account
that $-i\alpha(T_{\beta_j})\ge0$
for $\alpha\in\Delta^+_{\frak m}$
and $-i\beta_k(T_{\beta_j})=\delta^k_j$ (see subsection 4.1),
we complete the proof.
{}\hfill$\square$
\enddemo

Let
$\rho_{\frak k}$ be the mapping from
$(BC)_r^+\subset({\frak a}^{\Bbb C})^*$
to $({{\frak t}'}^{\Bbb C})^*$ with the graph given by
$$
\bigl\{
({\tfrac i2}\epsilon_j,
{\tfrac i2}\epsilon'_j),
(\epsilon_j,0),\ j=\overline{1,r};\
\bigl({\tfrac i2}(\epsilon_p\pm\epsilon_k),
{\tfrac i2}(\epsilon'_p-\epsilon'_k)\bigr),\ 1\le p<k\le r
\bigr\}.
$$

\proclaim{Corollary 4.7.3}
Let
$\lambda\in\Sigma^{+}$ be a restricted root. Suppose that
$\lambda_I\ne0$.
Let $\xi_\lambda\in{\frak m}_\lambda$ and
$\xi'_{\lambda_I}=I\xi_\lambda\in{\frak m}_{\lambda_I}$. Then
for any $T\in{\frak t}'$ and $\zeta_\lambda\in{\frak k}_\lambda$,
$\zeta'_{\lambda_I}\in{\frak k}_{\lambda_I}$ with the notations of
Lemma~4.3
$$
\operatorname{ad}_T(\zeta_\lambda)=
i\rho_{\frak k}(\lambda)(T)\zeta'_{\lambda_I}.
$$
If $\lambda_I=0$ then $[{\frak t}',{\frak k}_\lambda]=0$.
\endproclaim
\demo{Proof}
Using the notations of Lemma~4.3 and the Jacobi identity for
the vectors $X_j,IX_j$ and $\zeta_\lambda\in{\frak k}_\lambda$ we obtain
that
$$\align
[[X_j,IX_j],\zeta_\lambda]
&=[X_j,[IX_j,\zeta_\lambda]]-[IX_j,[X_j,\zeta_\lambda]]
=2[X_j,I[X_j,\zeta_\lambda]] \\
&=-2i\lambda(X_j)[X_j,I\xi_\lambda]
=2\lambda(X_j)\lambda_I(X_j)\zeta'_{\lambda_I}.
\endalign$$
If
$\lambda_I=0$, in the chain of equations above
$[X_j,I\xi_\lambda]=0$ because $I\xi_\lambda\in{\frak a}={\frak m}_0$.
Now the assertion of the corollary comes
from Proposition~4.7.
{}\hfill$\square$
\enddemo

By~(3.1)
$\operatorname{ad}_{Iw}^2(I\xi)=I\operatorname{ad}_w^2(\xi)$,
and by~(4.10)
$I\operatorname{ad}_{[w,Iw]}(\xi)
=(\operatorname{ad}_{w}^2+\operatorname{ad}_{Iw}^2)(\xi)$
(for all $\xi\in{\frak m}$). So as an immediate
consequence of Proposition~4.7 we obtain
\proclaim{Corollary 4.7.4}
Sum~{\rm(4.7)} {\rm(}resp.~{\rm(4.13))}
is an orthogonal eigenspace splitting of
${\frak m}$ for all operators $\operatorname{ad}_w^2$,
$\operatorname{ad}_{Iw}^2$, $w\in{\frak a}$
{\rm(}resp. $I\operatorname{ad}_{[w,Iw]}$,
$\operatorname{ad}_{w}^2+\operatorname{ad}_{Iw}^2$, $w\in{\frak a}${\rm)}.
In particular, for each
$w\in{\frak m}$,
$[\operatorname{ad}_w^2,\operatorname{ad}_{Iw}^2]|{\frak m}=0$ and
$I\operatorname{ad}_{[w,Iw]}
=\operatorname{ad}_{w}^2+\operatorname{ad}_{Iw}^2$ on ${\frak m}$.
\endproclaim

\remark{Remark 4.8}
The well-known
Harish-Chandra and Moore theorem
for Hermitian symmetric spaces (see~\cite{Wo})
describes the restricted root system of
$({\frak g},{\frak t}')$ and, using
the Cayley transform, such a system for
$({\frak g},{\frak a})$ (the spaces
${\frak t}'$ and
${\frak a}$ are conjugated in
${\frak g}$ as Cartan subalgebras of the same compact Lie algebra
${\frak a}\oplus I{\frak a}\oplus[{\frak a},I{\frak a}]$). This
theorem follows from Corollaries~4.7.2
and~4.7.3. But the mapping
$\Sigma^+\to\Sigma^+\cup\{0\}$,
$\lambda\mapsto\lambda_I$ of Proposition~4.7 allows us to find
direct connection between the restricted root decompositions for
$({\frak g},{\frak a})$ and $({\frak g},{\frak t}')$.
\endremark

{\bf 4.3. The main theorem.}
Here using the result of previous subsection we
construct all antiholomorphic $K$-equivariant mappings
on homogeneous domains in ${\frak m}$
and prove the main theorem.

Let
$\sigma:T(G/K)\to T(G/K)$ be the
involution which maps any tangent vector
$Y$ at $gK$ onto $-Y$ at $gK$ and let
$P\in \text{Alm}(W)$, where
$W=-W$. It is easy to see that
$\sigma$ is an antiholomorphic
involution for the almost complex structure
$J(P)$, i.e.
$\sigma_*(F(P))=\overline{F(P)}$,
iff $\overline{P_w}=P_{-w}$ for all
$w\in W$. But
$I\in\operatorname{Ad}(K)|{\frak m}$ (see subsection 4.1), so
by~(2.6) $IP_w=P_{Iw}I$ and, consequently,
$P_w=P_{-w}$. We have proved
\proclaim{Lemma 4.9}
Let $P\in \text{Alm}(W)$.
The mapping
$\sigma$ is an antiholomorphic involution for
$J(P)$ iff
$\overline P=P$.
\endproclaim

The following proposition and lemma will be crucial for the
subsequent part of the paper.
\proclaim{Proposition 4.10}
Let $P=(R+iS)\in\text{Alm}(W)$ and
$RI=IR$, $SI=-IS$. Suppose that
$L_\xi(SR^{-1})\eta=L_\eta(SR^{-1})\xi$ for all {\rm(}fixed{\rm)}
$\xi,\eta\in{\frak m}$. Then
$SR^{-1}=(a_1+a_2I)\Upsilon_*$ on $W$,
where $a_1,a_2\in{\Bbb R}$ are some numbers and
$\Upsilon$ is a rational $K$-equivariant mapping on
${\frak m}$ given by
$$
\Upsilon(w)=(-\operatorname{ad}_w^2
-\operatorname{ad}^2_{Iw})^{-1}(w).
$$
The endomorphism $\Upsilon_{*w}:{\frak m}\to{\frak m}$, where
$w$ belongs to the set $W_\Upsilon$  of all regular points of
$\Upsilon$, anticommute with
$I$ iff all restricted roots from
$\Sigma$ are indivisible, i.e.
$\Sigma$ has type
$C_r$.
\endproclaim
\demo{Proof}
Fix the Cartan subspace
${\frak a}=\sum_{j=1}^{r}{\Bbb R} X_{\beta_j}$
of ${\frak m}$. Let ${\frak A}={\frak a}\oplus I{\frak a}$.
Consider $({\frak m},I)$ as a space over
${\Bbb C}$ with fixed basis
$\{X_\alpha,\alpha\in\Delta^+_{\frak m}\}$.
For each complex number $z_0=x+\imath y$ put
$z_0X_\alpha\overset\text{def}\to{=}
xX_\alpha+yIX_\alpha =xX_\alpha +yY_\alpha\in{\frak m}$
(here $\imath^2=-1$,
$x,y\in{\Bbb R}$). Denote by
$C:{\frak m}\to{\frak m}$ the corresponding
complex conjugation mapping, i.e.
$C(z_0X_\alpha)=\overline{z_0}X_\alpha$, where
$\overline{z_0}\overset\text{def}\to{=} x-\imath y$.
For any complex vector
$z=(z_1,..,z_r)\in{\Bbb C}^r$ denote by
$zX$ the vector
$\sum_{j=1}^{r}z_jX_j\in{\frak A}$. Put
$Z=\{z\in{\Bbb C}^r: zX\in W\cap{\frak A}\}$. Since
$ICSR^{-1}=CSR^{-1}I$, there exists the complex matrix-function
$z\mapsto(h_{\alpha\beta}(z))$,
$\alpha,\beta\in\Delta^+_{\frak m}$ on
$Z$ such that the ${\Bbb C}$-linear mapping
$C(SR^{-1})_{zX}$ takes each
$\sum_\beta v_\beta X_\beta\in{\frak m}$ to
$\sum_{\alpha,\beta} h_{\alpha\beta}(z)v_\beta\cdot X_\alpha\in{\frak m}$.
Considering
$C$ as an ${\Bbb R}$-linear mapping on the set
$W\cap{\frak A}$ and using relation~(3.9),
we conclude that locally each function
$h_{\alpha\beta}:Z\to{\Bbb C}$ is a partial
derivative of some holomorphic function, i.e.
$h_{\alpha\beta}$ is holomorphic on the set
$Z\subset{\Bbb C}^r$. In particular, each holomorphic 1-form
$\sum_{j=1}^{r}h_{\alpha\beta_j}(z)\cdot dz_j$ is closed.
From~(4.1)
it follows that for any vector $T\in{\frak t}'$,
$\exp(\operatorname{ad}_{T})(X_\alpha)=e^{\imath \alpha'(T)}X_\alpha$,
where
$\alpha'=-i\alpha|{\frak t}'$ is a real linear function on
${\frak t}'$. Because of $K$-equivariance of $SR^{-1}$
$$
\operatorname{Ad}_k C(CSR^{-1})_w
\operatorname{Ad}_{k^{-1}}=
C(CSR^{-1})_{\operatorname{Ad}_k w},
\qquad \forall w\in W, k\in K,
$$
and, consequently,
$$
e^{-\imath(\alpha'+\beta')(T)}
h_{\alpha\beta}(z_1,..,z_r)=
h_{\alpha\beta}(e^{\imath\beta'_1(T)}z_1,..,
e^{\imath\beta'_r(T)}z_r),
\qquad\forall z\in Z,\ T\in{\frak t}'.
$$
Therefore
$h_{\alpha\beta}(z_1,..,z_r)
=a_{\alpha\beta}\cdot z_1^{p_1}\cdots z_r^{p_r}$
if $\alpha'+\beta'=-\sum_{k=1}^{r}p_k\beta'_k$
for some integers $p_k\in{\Bbb Z}$ else
$h_{\alpha\beta}(z)=0$. Here
$a_{\alpha\beta}\in{\Bbb C}$ is
some constant. Then by Corollary~4.7.2
$$
h_{\alpha\beta}(z)=\cases
a_{\alpha\beta}\cdot z_k^{-1}z_j^{-1}, &\ \text{if}\quad
\alpha=\beta_k,\ \beta=\beta_j,\
k,j=\overline{1,r};\\
a_{\alpha\beta}\cdot z_p^{-1}z_k^{-1}, &\ \text{if}\quad
\alpha'=\beta'=\frac12(\epsilon'_p+\epsilon'_k),
\ p\ne k;\\
a_{\alpha\beta}\cdot z_j^{-1}, &\ \text{if}\quad
\alpha'=\beta'=\frac12\epsilon'_j,\ j=\overline{1,r};
\endcases
$$
and otherwise $h_{\alpha\beta}(z)=0$.

Let
${\frak m}_{1}$ be a subspace of
${\frak m}$ spanned by vectors
$X_\alpha,Y_\alpha$ with
$\alpha'=\frac12(\epsilon'_p+\epsilon'_k)$,
$1\le p\le k\le r$. Denote by ${\frak m}_{1/2}$
the orthogonal complement to the space ${\frak m}_{1}$
in ${\frak m}$. Using the expressions for the
matrix element $h_{\alpha\beta}$ of $CSR^{-1}$, we
obtain that the subspaces ${\frak m}_{1}$ and
${\frak m}_{1/2}$ are invariant with respect
to all operators $(CSR^{-1})_w$ with
$w\in W\cap{\frak a}$ and
$$\align
\tfrac{d}{dt}\bigr|_{0}(CSR^{-1})_{w+tw}|{\frak m}_{1}
&=-2(CSR^{-1})_w|{\frak m}_{1}; \\
\tfrac{d}{dt}\bigr|_{0}(CSR^{-1})_{w+tw}|{\frak m}_{1/2}
&=-(CSR^{-1})_w|{\frak m}_{1/2}.
\endalign$$
Consider the mapping
$$
{\tilde\Upsilon}: W\to{\frak m},\quad
{\tilde\Upsilon}(w)=-(SR^{-1})_w\cdot w
$$
This mapping is
$K$-equivariant because the mapping
$w\mapsto (SR^{-1})_w$ satisfies condition~(2.6). By~(3.9) for
$\xi\in{\frak m}$
$$
{\tilde\Upsilon}_{*w}(\xi)=-(L_\xi(SR^{-1}))_w w-(SR^{-1})_w\xi
=-(\tfrac{d}{dt}\bigr|_{0}(SR^{-1})_{w+tw})\xi-(SR^{-1})_w\xi.
$$
Thus for all $w\in W\cap{\frak a}$
$$
{\tilde\Upsilon}_{*w}|{\frak m}_{1}
=(SR^{-1})_w|{\frak m}_{1}
\quad\text{and}\quad
{\tilde\Upsilon}_{*w}|{\frak m}_{1/2}=0.
\tag 4.17
$$

Let us find this function
${\tilde\Upsilon}$. Since the holomorphic form
$\sum_{j=1}^{r}h_{\beta_k\beta_j}\cdot dz_j$ is closed,
$h_{\beta_k\beta_j}=0$ if $k\ne j$.
Relabel the function $h_{\beta_j\beta_j }$ by $h_j$.
Since
$G/K$ is an irreducible hermitian symmetric
space, the restricted root system of
$({\frak g},{\frak a})$ has the type $C_r$ or
$(BC)_r$~\cite{12, Ch.X, \S 6}.
Therefore the restricted Weyl group of
$({\frak g},{\frak a})$ induces all signed permutations
$X_j\mapsto \pm X_{k(j)}$.
Taking into account
$K$-equivariance~(2.6)
again, we obtain that
$h_{s(j)}(z_1,..,z_r)=h_j(z_{s(1)},..,z_{s(r)})$,
where $s$ is any permutation. But the function
$h_j(z)$ depends only on the $j$-th coordinate $z_j$ of $z$.
So there is a unique function $h(t)=z_0t^{-2}$, $z_0,t\in{\Bbb C}$
such that $h_j(z)=h(z_j)$. Then
$(SR^{-1})_{zX}$ takes
each $v_jX_j\in{\frak A}$ to
$\overline{z_0z_j^{-2}v_j}\cdot X_j\in{\frak A}$,
i.e. ${\tilde\Upsilon}(\sum_j x_jX_j)=(a_1+a_2I)\sum_j x_j^{-1}X_j$,
where $-\overline{z_0}=a_1+\imath a_2$, $x_j\in{\Bbb R}$.

The restriction of the mapping
$w\mapsto(-\operatorname{ad}_w^2-\operatorname{ad}^2_{Iw})^{-1}(w)$ to
${\frak a}$ takes each $\sum_j x_jX_j$ to
$\sum_j x_j^{-1}X_j$. Thus the mappings
${\tilde\Upsilon}$ and $(a_1+a_2I)\Upsilon$ coincide on the set $W\cap{\frak a}$
and because of $K$-equivariance of ${\tilde\Upsilon}$ and $\Upsilon$
$$
{\tilde\Upsilon}(w)=(a_1+a_2I)\Upsilon(w)
\quad\text{for all}\quad w\in W.
\tag 4.18
$$

Splitting~(4.7) is the common eigenspace splitting for all
$\Upsilon_{*w}$,
$w\in{\frak a}$ (see Proposition~4.6). The restrictions of
$\Upsilon_{*w}$ and $I$ to
${\frak A}={\frak a}\oplus I{\frak a}$ anticommute because
by~(4.4) and~(4.9)
$$
\Upsilon_{*w}(X_j)=-x_j^{-2}X_j
\qquad\text{and}\qquad
\Upsilon_{*w}(IX_j)=x_j^{-2}IX_j,
\qquad w\in W\cap{\frak a},\ j=\overline{1,r}.
$$
Suppose that
$\lambda_{pk}^\pm=i(\epsilon_p\pm\epsilon_k)/2$,
$p\ne k$, $\lambda_{j/2}=i\epsilon_j/2$,
$1\le j\le r$ are restricted roots and
$\xi_{pk}^\pm\in{\frak m}_{\lambda_{pk}^\pm}$,
$\xi_{j/2}\in{\frak m}_{\lambda_{j/2}}$. Then applying~(4.4) to the
$K$-equivariant mapping
$\Upsilon$, we obtain that
$$
\Upsilon_{*w}(\xi_{pk}^\pm)
=\frac{x_p^{-1}\pm x_k^{-1}}{x_p\pm x_k}\xi_{pk}^\pm
=\pm\frac1{x_px_k}\xi_{pk}^\pm
\quad\text{and}\quad
\Upsilon_{*w}(\xi_{j/2})
=\frac1{x_j^2}\xi_{j/2},
\tag 4.19$$
where $w=\sum_{j=1}^r x_jX_j$.
But ${\tilde\Upsilon}_{*w}|{\frak m}_{1/2}=0$ (see~(4.17)).
Therefore by~(4.18) and~(4.19)
$a_1=a_2=0$ if ${\frak m}_{1/2}\ne0$, i.e.
$S=0$. Thus always ${\tilde\Upsilon}_*=SR^{-1}$.
Taking into account that
$I({\frak m}_{\lambda_{pk}^\pm})={\frak m}_{\lambda_{pk}^\mp}$ and
$I({\frak m}_{\lambda_{j/2}})={\frak m}_{\lambda_{j/2}}$,
we complete the proof.
{}\hfill$\square$
\enddemo

We can supplement Proposition~4.10 with the following
simple statement.
\proclaim{Corollary 4.10.1}
Suppose that the restricted root system
$\Sigma$ of
$({\frak g},{\frak a})$ has type
$C_r$. Then for arbitrary
$a_1,a_2\in{\Bbb R}$ the mapping
$(a_1+a_2I)\circ\Upsilon$ satisfies~{\rm(4.5)}.
\endproclaim
\demo{Proof}
Since
$\Upsilon\in\text{EC}(W_\Upsilon)$, it only
remains to prove~(4.5) for the mapping
$I\circ\Upsilon$. But
$I\Upsilon_{*w}=\Upsilon_{*Iw}I$ because the mapping
$\Upsilon$ is
$K$-equivariant and
$I$ and the automorphism
$\exp\frac\pi2\operatorname{ad}_{Z_0}\in\operatorname{Ad}(K)$
coincide when restricted to ${\frak m}$
(see subsection
4.1). By Proposition~4.10
$I\Upsilon_{*w}=-\Upsilon_{*w}I$. So that
$\Upsilon_{*Iw}=-\Upsilon_{*w}$. Taking into
account already proved identity~(4.5) for
$\Upsilon$ at
$Iw$ and properties~(3.1) of
$I$, we obtain
$$\align
I\Upsilon_{*w}([w,[\xi,\eta]])
&=\Upsilon_{*Iw}([Iw,[\xi,\eta]]) \\
&=[[Iw,\Upsilon_{*Iw}(\xi)],\eta]-
[[Iw,\Upsilon_{*Iw}(\eta)],\xi] \\
&=[[w,-I\Upsilon_{*Iw}(\xi)],\eta]-
[[w,-I\Upsilon_{*Iw}(\eta)],\xi] \\
&= [[w,I\Upsilon_{*w}(\xi)],\eta]-
[[w,I\Upsilon_{*w}(\eta)],\xi].
\endalign$$
\enddemo

\proclaim{Lemma 4.11}
Let
$(J(P),\Omega)$ be a K\"ahler structure on
$D$ such that $IP=\overline PI$ on
$W$. Then for each Cartan subspace
${\frak a}=\sum_{j=1}^r{\Bbb R} X_j\subset{\frak m}$,
$R_w({\frak a})={\frak a}$, where
$w\in W\cap{\frak a}$. Moreover,
$R_w(X_j)\in{\Bbb R} X_j$,
$j=\overline{1,r}$ and
$R_w({\frak M}_\lambda)\subset{\frak M}_\lambda$
for all roots $\lambda\in\Sigma^{++}$.
\endproclaim
\demo{Proof}
For arbitrary
mapping
$A\in\text{Eqv}(W)$ it follows from~(2.6) that
$[\operatorname{ad}_\zeta,A_w]=(L_{[\zeta,w]}A)_w$, where
$\zeta\in{\frak k}$. Then by~(3.8)
$$
\bigl[\operatorname{ad}_\zeta,(1+(SR^{-1})^2_w)R_w\bigr](\eta)
=I\Bigl[\eta,\bigl[w,R^{-1}_wI[\zeta,w]\bigr]\Bigr],
\qquad \forall \zeta\in{\frak k},\eta\in{\frak m}.
\tag 4.20$$

We first shall prove that
$R_w({\frak a})\subset{\frak A}$, where
${\frak A}={\frak a}\oplus I{\frak a}$. To see this, denote by
$K^{\frak a}$ the connected subgroup of
$K$ with the Lie algebra
${\frak k}^{\frak a}$ (the centralizer of the Cartan subspace
${\frak a}$ in
${\frak k}$). All automorphisms
$\operatorname{Ad}_k$,
$k\in K^{\frak a}$ leave the space
${\frak a}\oplus I{\frak a}$ pointwise fixed.
Then by~(2.6) all these automorphisms leave
$R_w({\frak a})$ pointwise fixed.
Now it follows from Proposition~4.2 that
$R_w({\frak a})\subset{\frak A}$ if
${\frak k}^{\frak a}\ne0$.

It remains to consider the case when
${\frak k}^{\frak a}=0$, i.e.
${\frak a}\subset{\frak m}$ is a Cartan subalgebra of
${\frak g}$ and
$\Sigma$ is a root system of
$({\frak g},{\frak a})$. Then
${\frak t}'={\frak t}$ and by Corollary~4.7.2 for each root
$\alpha\in\Delta^+_{\frak m}$ there is a unique root
$\lambda^+\in\Sigma^{++}$ such that
${\Bbb R} X_\alpha\oplus {\Bbb R} Y_\alpha=
{\frak m}_{\lambda^+}\oplus I{\frak m}_{\lambda^+}$.
Here
$\alpha=\frac i2(\epsilon'_p+\epsilon'_k)$ and
$\lambda^+=\frac i2(\epsilon_p+\epsilon_k)$, where
$1\le p\le k\le r$. If
$\lambda^+\ne i\epsilon_p$ then
$I{\frak m}_{\lambda^+}={\frak m}_{\lambda^-}$, where
$\lambda^-=\frac i2(\epsilon_p-\epsilon_k)$, else
$I{\frak m}_{\lambda^+}={\Bbb R} X_p$.
In this one-dimensional subspace
$I{\frak m}_{\lambda^+}$ fix a non-zero vector
$X_\alpha^-$ assuming that $X_{\beta_j}^-=X_j$,
$j=\overline{1,r}$. Then for
$\alpha=\frac i2(\epsilon'_p+\epsilon'_k)$,
$T=\sum_{j=1}^r t_jT_j\in{\frak t}'$,
$w=\sum_{j=1}^r x_jX_j\in{\frak a}$
$$\align
\operatorname{ad}_T(X_\alpha^-)
&={\tfrac 12}
(t_p+t_k)IX_\alpha^-,\\
\operatorname{ad}_w^2(X_\alpha^-)
&={\textstyle-\frac 14}
(x_p-x_k)^2X_\alpha^-,\tag 4.21 \\
\operatorname{ad}_w^2(IX_\alpha^-)
&={\textstyle-\frac 14}
(x_p+x_k)^2IX_\alpha^-.
\endalign$$
Now consider $({\frak m},I,\langle ,\rangle)$ as an Hermitian
space over ${\Bbb C}$ with fixed orthogonal basis
$\{X_\beta^-,\beta\in\Delta^+_{\frak m}\}$. Assume
that all base vectors have the same length which is equal
to length of (conjugated) vectors $\{X_j\}_{j=1}^r$.
Since each operator
$R_w:{\frak m}\to{\frak m}$ is symmetric and commuting with
$I$, the corresponding complex matrix is Hermitian. Let
$(R_{\beta j})$ and $(r_{\beta j})$,
$\beta\in\Delta^+$,
$j=\overline{1,r}$ be two complex matrix-functions
corresponding to the operator-functions
$R|{\frak A}$ and $R^{-1}|{\frak A}$. Then
$\sum_{\beta\in\Delta^+_{\frak m}}\overline{R_{\beta p}}
r_{\beta k}=\delta^p_k$.

Put $\eta=X_n\in{\frak a}$ and
$\zeta=T\in{\frak t}'$ in~(4.20). By Proposition~4.10
$(SR^{-1})^2_w=(a_1^2+a_2^2)(\Upsilon_*)_w^2$. Taking into account
relations~(4.21) and~(4.19), we obtain from~(4.20) the following
equation
$$\align
&
\bigl({\tfrac12} t_p+{\tfrac12} t_k-t_n\bigr)
\bigl(1+(a_1^2+a_2^2)x_p^{-2}x_k^{-2}\bigr)R_{\alpha n}
\tag 4.22 \\
&=\sum_{j=1}^r t_jx_j\bigl(
{\tfrac14}(\delta_p^n-\delta_k^n)(x_p-x_k)\hat r_{\alpha j}
+{\tfrac\imath 4}(\delta_p^n+\delta_k^n)(x_p+x_k)
\check r_{\alpha j}\bigl)
\endalign$$
Here
$\alpha=\frac i2(\epsilon'_p+\epsilon'_k)$,
$\hat r_{\alpha j}$ and
$\check r_{\alpha j}$ are real and
imaginary parts of the complex function
$r_{\alpha j}$. This equation then gives that
$R_{\alpha n}=r_{\alpha n}=0$ if
$n\ne p$ and
$n\ne k$.

Assume now that
$p\ne k$. It follows from~(4.22) that
$R_{\alpha p}=\overline{R_{\alpha k}}$ and
$$\align
R_{\alpha p}\bigl(1+(a_1^2+a_2^2)x_p^{-2}x_k^{-2}\bigr)
&=-{\tfrac12}x_p\bigl((x_p-x_k)\hat r_{\alpha p}
+\imath(x_p+x_k)\check r_{\alpha p}\bigr)\tag 4.23\\
&={\tfrac12}x_k\bigl((x_p-x_k)\hat r_{\alpha k}
+\imath(x_p+x_k)\check r_{\alpha k}\bigr).
\endalign$$
But
$$
\sum_{\beta\in\Delta^+_{\frak m}}\overline{R_{\beta p}}r_{\beta k}
=\overline{R_{\alpha p}}r_{\alpha k}=\delta_k^p=0,
$$
i.e. either $R_{\alpha p}$ or
$r_{\alpha k}$ equal zero.
By~(4.23) these two functions on
$W\cap{\frak a}$ equal zero simultaneously. Thus
$R_w({\frak a})\subset{\frak A}$.

Turning to the general case, consider
again relation~(4.20). Using the basis
$\{X_j\}$, $j=\overline{1,r}$ in
${\frak a}$, we deduce the following
equations for the matrix elements
$R_{pn}$ and
$r_{pn}$ of
$R,R^{-1}:{\frak a}\to{\frak A}$
$$
(t_p-t_n)\bigl(1+(a_1^2+a_2^2)x_p^{-4}\bigr)R_{pn}
=\imath\delta_p^n{\textstyle\sum_{j=1}^r}
t_jx_jx_p(\operatorname{Im}r_{pj}),
\qquad \forall t_j\in{\Bbb R},\ j=\overline{1,r}.
$$
Thus $R_{pn}=0$ if $n\ne p$ and
$r_{pp}\in{\Bbb R}$, i.e.
$R_w({\Bbb R} X_p)\subset{\Bbb R} X_p$. But
$I[{\frak t}',{\frak a}]\subset{\frak a}$. Therefore for any
$\zeta\in{\frak t}'$ the expression on the
right in (4.20) vanishes and, consequently,
$(1+(SR^{-1})^2_w)R_w({\frak M}_\lambda)\subset{\frak M}_\lambda$.
{}\hfill$\square$
\enddemo

We showed above that a structure of the mapping
$P$ depends on a type of the restricted root system
$\Sigma$ of the symmetric space
$G/K$. For each type we define a set
${\Cal A}^\Sigma\subset{\Bbb R}^3\times\{\pm1\}$ by
$$
{\Cal A}^{C}={\Bbb R}^3\times\{1\}
\quad\text{and}\quad
{\Cal A}^{BC}={\Bbb R}^+\times\{0\}\times\{0\}
\times\{\pm1\}\cup\{(0,0,0,-1)\}.
\tag 4.24$$
For an element ${\bold a}\in{\Cal A}^{C}$, let
$a_\dagger=\frac12\bigl(\sqrt{a_0^2+4a_1^2+4a_2^2}-a_0\bigr)$
if $a_1^2+a_2^2>0$ and $a_\dagger=-a_0$ if
$a_1=a_2=0$. Put
${\frak a}_\dagger=\{\sum_{j=1}^r
x_j X_j\in{\frak a}: x_j^2>a_\dagger\}$.
This subset of the Cartan subspace ${\frak a}$ defines a unique
$\operatorname{Ad}(K)$-invariant open connected subset
$W_{\bold a}^{C}$ of ${\frak m}$ with
$W_{\bold a}^{C}\cap{\frak a}={\frak a}_\dagger$.
For ${\bold a}\in{\Cal A}^{BC}$, let
$W_{\bold a}^{BC}={\frak m}$ if $\varepsilon=1$
and $W_{\bold a}^{BC}={\frak m}\setminus\{0\}$
if $\varepsilon=-1$. Set
$D_{\bold a}^\Sigma=\Pi(G\times W_{\bold a}^\Sigma)$.

The central result of this paper is the following.
\proclaim{Theorem 4.12}
Let
$(J(P),\Omega)$ be a $G$-invariant K\"ahler structure on the
$G$-invariant domain
$D\subset T(G/K)$, where
$G/K$ is the irreducible Hermitian symmetric
space of compact type. Suppose that
$IP=\overline PI$ on
$W$. Then there exists a unique quadruple
$(a_0,a_1,a_2,\varepsilon)\in{\Cal A}^\Sigma$ such that
$W\subset W_{\bold a}^{\Sigma}$ and
$$
P_w=\bigl(1+i(a_1{+}a_2I)\Upsilon_{*w}\bigr)\cdot
\bigl(1+(a_1^2{+}a_2^2)(\Upsilon_*)_w^2\bigr)^{-1}\cdot B_w,
\quad w\in W\subset{\frak m},
\tag 4.25$$
$$
B_w=\bigl(I\cdot\operatorname{ad}_{[I\hat\phi(w),w]}
+\varepsilon\sqrt{|a_0|}\cdot \operatorname{Id}\bigr)|{\frak m},
\quad \hat\phi(w)=\phi(-\operatorname{ad}_w^2
-\operatorname{ad}_{Iw}^2)w,
\tag 4.26$$
$$
\phi(t)=\frac{\sqrt{a_0+t-(a_1^2+a_2^2)t^{-1}}-
\varepsilon\sqrt{|a_0|}}{t}.
$$

For arbitrary
$(a_0,a_1,a_2,\varepsilon)\in{\Cal A}^\Sigma$ the operator-function
$P$~{\rm(4.25)} determines a K\"ahler structure
$(J(P),\Omega)$ on the $G$-invariant domain
$D_{\bold a}^{\Sigma}\subset T(G/K)$. This structure
anticommutes with $J^-$.

Moreover, if $a_2=0$, this K\"ahler structure
$(J(P),\Omega=d\theta)$~{\rm(4.25)} admits a potential function $Q$, i.e.
$\Omega=2i\overline\partial\partial Q$; if, in addition,
$a_1=0$, then $\theta=2\operatorname{Im}\overline\partial Q$.
The function  $(\Pi^*Q)(g,w)=\langle q(-\operatorname{ad}_w^2
-\operatorname{ad}_{Iw}^2)w,w\rangle$,
where
${\dsize q(t)=\frac{1}{2t}\int\frac{dt}{\sqrt{a_0+t-a_1^2t^{-1}}}}$.
\endproclaim
\demo{Proof}
By Proposition~4.10
$SR^{-1}=(a_1+a_2I)\Upsilon_*$ on
$W$. Define the $K$-equivariant mapping
$B\in\text{Eqv}(W)$ putting
$B_w=(1+c^2(\Upsilon_*)_w^2) R_w$,
$c^2=a_1^2+a_2^2$. A change
$R_w\mapsto (1+c^2(\Upsilon_*)_w^2)^{-1}B_w$
converts~(3.8) into
$$
(L_\xi B)_w(\eta)
=-[[w,IB_w^{-1}(1+c^2(\Upsilon_*)_w^2)\xi],I\eta],
\quad w\in W,\ \xi,\eta\in{\frak m}.
\tag 4.27$$
By Lemma~4.11
$B_w(X_j)=b_j(w)X_j$ for all
$w\in W\cap{\frak a}$. Putting in~(4.27)
$\eta=X_j$ and
$\xi=X_k$, we obtain that the function
$b_j$ on $W\cap{\frak a}$ depends only on the
$j$-th coordinate of the vector
$w=\sum_{j=1}^r x_jX_j$. Taking into account
the action of the restricted Weyl group of
$({\frak g},{\frak a})$ on ${\frak a}$ and
$K$-equivariance of
$B$, we conclude that all
$b_j$ coincide as functions on some subset of
${\Bbb R}$ (see the proof of Proposition~4.10
and~(2.6)). This unique function will be denoted by
$b$, i.e.
$b_j(w)=b(x_j)$.
Solving equation~(4.27) for $\xi=\eta=X_j$, i.e.
$b'(x)=b^{-1}(x)(1+c^2x^{-4})x$, we find
$$
b(x)=\sqrt{a_0+x^2-(a_1^2+a_2^2)x^{-2}}
\tag 4.28$$
for some constant $a_0\in{\Bbb R}$.
Also $B_w(IX_j)=b_j(w)IX_j$ because $B_wI=IB_w$.

By Lemma~4.11 $B_w({\frak M}_\lambda)={\frak M}_\lambda$ for each
$\lambda\in\Sigma^{++}$. Since
$[{\frak a},I{\frak a}]\subset{\frak t}'$, from~(4.27) and
Corollary~4.7.2 it follows that~(4.13) is the
orthogonal eigenspace splitting for all operators
$(L_\xi B)_w$ with $w\in W\cap{\frak a}$,
$\xi\in{\frak a}$, and, consequently,
there exists a constant operator
$B_{\dagger}:{\frak m}\to{\frak m}$
such that
$(B_w-B_{\dagger})|{\frak M}_\lambda=b_\lambda(w)\cdot
\operatorname{Id}_{{\frak M}_\lambda}$. We can choose
the operator $B_{\dagger}$ with trace $=0$ when restricted
to each ${\frak M}_\lambda$.

Put $B_w^\lambda=B_w|{\frak M}_\lambda$ for
$w\in W\cap{\frak a}$.
Fix an element $\lambda\in\Sigma^{++}$ and a vector
$\xi_\lambda\in{\frak m}_\lambda$. There is a vector
$\zeta_\lambda\in{\frak k}_\lambda$ satisfying~(4.3).
The change $R_w\mapsto (1+c^2(\Upsilon_*)_w^2)^{-1}B_w$
converts~(4.20) into
$$
\bigl[\operatorname{ad}_\zeta,B_w\bigr](\eta)
=I\Bigl[\eta,\bigl[w,B_w^{-1}
(1+c^2(\Upsilon_*)_w^2)I[\zeta,w]\bigr]\Bigr],
\qquad \forall \zeta\in{\frak k},\eta\in{\frak m}.
\tag 4.29$$
We claim that $B_w({\frak m}_\lambda)\subset{\frak m}_\lambda$.
Indeed, if
$\lambda=\lambda_I$ then ${\frak M}_\lambda={\frak m}_\lambda$,
i.e. the claim holds.
Suppose now that $\lambda\ne\lambda_I$ and $\lambda_I\ne0$. Then
${\frak M}_\lambda={\frak m}_\lambda\oplus{\frak m}_{\lambda_I}$.
Put in~(4.29)
$\zeta=\zeta_\lambda$ and
$\eta=w'\in{\frak a}$.
In view of~(4.19)
$\Upsilon_{*w}(\xi_\lambda)=u_\lambda(w)\xi_\lambda$.
Because of Lemma~4.3
equation~(4.29) gives
$$\align
\lambda(B_ww')\xi_\lambda
&-\lambda(w')B_w(\xi_\lambda)
=-\lambda(w)\lambda_I(w)\lambda_I(w')
(1+c^2 u_\lambda^2(w))\cdot
(B_w^{-1}\xi_\lambda)_{{\frak m}_\lambda} \\
&-\lambda(w)\lambda(w)\lambda(w')
(1+c^2 u_\lambda^2(w))\cdot
(B_w^{-1}\xi_\lambda)_{{\frak m}_{\lambda_I}}.
\tag 4.30
\endalign$$
Setting equal
in~(4.30) the components belonging to the subspace
${\frak m}_{\lambda_I}$, we obtain that
$$
(B_w\xi_\lambda)_{{\frak m}_{\lambda_I}}
=\lambda^2(w)(1+c^2 u_\lambda^2(w))\cdot
(B_w^{-1}\xi_\lambda)_{{\frak m}_{\lambda_I}}.
$$
Similarly for the vector
$\xi_{\lambda_I}=I\xi_\lambda\in{\frak m}_{\lambda_I}$
we have
$$
(B_w\xi_{\lambda_I})_{{\frak m}_\lambda}
=\lambda^2_I(w)(1+c^2 u_{\lambda_I}^2(w))\cdot
(B_w^{-1}\xi_{\lambda_I})_{{\frak m}_\lambda}.
$$
But
$(B_wI\xi_\lambda)_{{\frak m}_\lambda}=
I(B_w\xi_\lambda)_{{\frak m}_{\lambda_I}}$
and an analogous relation holds for the operator $B_w^{-1}$.
Note also that  $u_\lambda^2=u_{\lambda_I}^2$ and
$\lambda^2\ne\lambda_I^2$. Therefore
$(B_w\xi_\lambda)_{{\frak m}_{\lambda_I}}=0$. The
claim is proved.

Turning to the general case,
we obtain for arbitrary
$\lambda\in\Sigma^{++}$ that
$$
\lambda(B_ww')\xi_\lambda
-\lambda(w')B_w(\xi_\lambda)
=-\lambda(w)\lambda_I(w)\lambda_I(w')
(1+c^2 u_\lambda^2(w))
B_w^{-1}(\xi_\lambda).
$$
Moreover, we can replace the
vector $\xi_\lambda$ in this equation by the vector $I\xi_\lambda$
because $B_wI=IB_w$. Thus
$$
\lambda(B_ww')\cdot B_w^\lambda
-\lambda(w')\cdot (B_w^\lambda)^2
=-\lambda(w)\lambda_I(w)\lambda_I(w')
(1+c^2 u_\lambda^2(w))\cdot
\operatorname{Id}_{{\frak M}_\lambda}.
$$
Since
$B_ww'=\sum_{j=1}^r b(x_j)x'_jX_j$,
the linear functionals
$\lambda(B_ww')$ and $\lambda(w')$ on
${\frak a}$ are linearly independent
for any $w$ in general position. Thus $B_\dagger=0$ and
$B_w^\lambda$ is a scalar operator. Then the real function
$b_\lambda(w)$ satisfies the equation
$$
\lambda(B_ww')\cdot b_\lambda(w)
-\lambda(w')\cdot b_\lambda^2(w)
=-\lambda(w)\lambda_I(w)\lambda_I(w')
(1+c^2 u_\lambda^2(w)).
$$
This equation has the following solutions
$$
b_\lambda={\tfrac12}(b_p+b_k)
\ \text{if}\ \lambda={\tfrac i2}(\epsilon_p+\epsilon_k)
\ \text{and}\
b_\lambda={\tfrac12}(b_j+\varepsilon\sqrt{a_0})
\ \text{if}\ \lambda={\tfrac i2}\epsilon_j,
\tag 4.31$$
where $\varepsilon=\pm1$. In other words,
with restrictions~(4.24) on the parameters
$a_0,a_1,a_2,\varepsilon$
$$
B_w=(I\operatorname{ad}_{T(w)}+\varepsilon\sqrt{|a_0|}
\operatorname{Id})|{\frak m},
\ \text{where}\
T(w)=\sum_{j=1}^r-(b_j(w){-}\varepsilon\sqrt{|a_0|})T_j.
\tag 4.32$$
Since $[IX_j,X_j]=-T_j$ and
$-\operatorname{ad}_{Iw}^2(X_j)=x_j^2X_j$,
expressions~(4.26) and~(4.32) for $B_w$,
$w\in W\cap{\frak a}$ coincide. The equivariance of the mapping
$B:w\mapsto B_w$~(4.26) proves the first assertion of the theorem.

Since
$\Upsilon_{*w}I=-I\Upsilon_{*w}$, each operator
$1+i(a_1+a_2I)\Upsilon_{*w}$ is invertible.
So by construction the operator-function
$P$~(4.25) determines an almost complex structure on
$D_{\bold a}^{\Sigma}$ anticommuting with
$J^-$. Therefore to prove the second assertion of the theorem it
suffices to show that the almost complex structure
$J(P)$ with
$P$~(4.25) is integrable. The following lemma generalizes some
assertions of Lemma~6 of~\cite{BG1}, where the case
$\operatorname{Im} P=0$ was considered.
\proclaim{Lemma 4.13}
Let $P\in\text{Alm}(W)$, $RS=SR$ and $IP=\overline PI$. Suppose
that $P$ satisfies conditions~{\rm(3.8)} and~{\rm(3.9)}.
Then the almost complex structure $J(P)$ on $D$ is integrable.
\endproclaim
\demo{Proof}
The lemma follows immediately from Lemma~3.6, Proposition~4.10
and its Corollary~4.10.1.
{}\hfill$\square$
\enddemo

From Proposition~4.10 and the proof above it follows
that~(4.7) is eigenspace splitting for all operators
$(SR^{-1})_w$ and $B_w$,
$w\in W\cap{\frak a}$ and, consequently, for all
$S_w$ and $R_w$. By equivariance
$R_w,S_w$ is a pair of commuting operators for all
$w\in W$.

Thus to prove integrability of $J(P)$, we have to verify only
condition~(4.27). Without loss of generality, we may
assume that $w\in W\cap{\frak a}$. By equivariance of $B$,
this condition
is equivalent to~(4.27) with $\xi\in{\frak a}$ and~(4.29).

Considering equation~(4.27) with
$\eta=\xi_\lambda$,
$\lambda\in\Sigma^{++}$ and the restriction
$b_\lambda|{\frak a}$ as a function of
$x_1,..,x_r$, we obtain
$$
(\partial b_\lambda/\partial x_j)\cdot\xi_\lambda=
-I\operatorname{ad}_{T_j}(\xi_\lambda)
\cdot(1+c^2x_j^{-4})x_jb_j^{-1},
\quad j=\overline{1,r}.
\tag 4.33$$
By Corollary~4.7.2
$I\operatorname{ad}_{T_j}(\xi_\lambda)
=i\rho_{\frak m}(\lambda)(T_j)\xi_\lambda$,
so~(4.33) is a linear combination of
equations for $b_j$ with the solutions $b(x_j)$~(4.28).

Since $I$ commutes with $\operatorname{ad}_\zeta$, from~(4.32)
it follows that the left-hand side
of~(4.29) equals $I\operatorname{ad}_{[\zeta,T(w)]}$.
But the $\operatorname{ad}$-representation of ${\frak k}$
in ${\frak m}$ is faithful (and irreducible)~\cite{GG, (8.5.1)}.
Hence~(4.29) is equivalent to the equation
$$
[T(w),\zeta]
=\bigl[w,B_w^{-1}
(1+c^2(\Upsilon_*)_w^2)I[\zeta,w]\bigr],
\qquad \forall \zeta\in{\frak k}.
\tag 4.34$$
For $\zeta\in{\frak k}_0={\frak k}^{\frak a}$
the left and right sides of~(4.34) vanish, because
$w\in{\frak a}$ and
$T(w)\in{\frak t}'=[{\frak a},I{\frak a}]$.
Similarly we have zeroes for
$\zeta\in{\frak k}_\lambda$,
$\lambda\in\Sigma^+$, where
$\lambda_I=0$, because
$[{\frak t}',{\frak k}_\lambda]=0$ by Corollary~4.7.3 and
$I[{\frak a},{\frak k}_\lambda]\subset I{\frak m}_\lambda\subset{\frak a}$.

Suppose now that
$\zeta=\zeta_\lambda\in{\frak k}_\lambda$ and
$\lambda_I\ne0$. Applying Corollary~4.7.3 again we obtain that
$$
i\rho_{\frak k}(\lambda)(T(w))\cdot \zeta'_{\lambda_I}
=-\lambda(w)\lambda_I(w)(1+c^2u_\lambda^2)b_\lambda^{-1}
\cdot \zeta'_{\lambda_I}.
$$
Here
$\lambda=\frac i2(\epsilon_p\pm\epsilon_k)$, $p\ne k$ or
$\lambda=\frac i2\epsilon_j$. It is easy to verify this algebraic
identity using~(4.19), (4.28), (4.31) and the expressions for
$\rho_{\frak k}(\lambda)$~(4.14). Thus the almost complex structure
$J(P)$ is integrable.

To prove the last assertion of the theorem for
the complex structure
$(J(P),\Omega)$
consider its subbundle $F(P)$ of $(0,1)$-vectors.
By definition $\overline\partial Q|F(P)=dQ|F(P)$ and
$\overline\partial Q|\overline{F(P)}=0$.
Denote by $\Phi$ the one-form
$\Pi^*(\overline\partial Q)$ on $G\times W$.
Then for any
$\xi\in{\frak m}$, $\zeta\in{\frak k}$, $w\in{\frak m}$:
$$
\Phi_{(g,w)}(\xi^l(g),-i\overline
{P_w}(\xi))=0,
\qquad
\Phi_{(g,w)}(\zeta^l(g),[w,\zeta])=0
$$
and
$$
\Phi_{(g,w)}(\xi^l(g),iP_w(\xi))
= i\langle \hat q(w),P_w(\xi)\rangle+
i\langle \hat q_{*w}\bigl(P_w(\xi)\bigr),w\rangle.
$$
Fix $w\in W\cap{\frak a}$. Then
$\hat q(w)\in{\frak a}$,
$\hat q_{*w}({\frak a})\subset{\frak a}$ and
$P_w({\frak a})\subset{\frak a}^{\Bbb C}$ because
$a_2=0$. Now using the invariance of the space
${\frak a}^\bot\subset{\frak m}$ with respect to
$P_{w}$ and $\hat q_{*w}$, we obtain that
$\Phi_{(e,w)}({\frak k}\oplus{\frak a}^\bot,{\frak a}^\bot)=0$
($[w,{\frak k}]\bot{\frak a}$).
Since the endomorphisms
$P_w$ and
$\hat q_{*w}$ commute and are symmetric, we have
for $\xi_0\in{\frak a}$
$$
\Phi_{(e,w)}(\xi_0,-i(R-iS)_w\xi_0)=0
$$
and
$$
\Phi_{(e,w)}(\xi_0,i(R+iS)_w\xi_0)=
i\langle (R+iS)_w(\hat q(w)+\hat q_{*w}(w)),\xi_0\rangle.
$$
Thus for any
$\eta\in{\frak g}, u\in{\frak m}$ we have
$$
\Phi_{(e,w)}(\eta,u)={\tfrac i2}
\langle (1+SR^{-1}SR^{-1})_wR_w
(w'),\eta\rangle+
{\tfrac 12}\langle (1+iSR^{-1})_w
(w'), u\rangle.
\tag 4.35$$
Here $w'$ is the vector
$\hat q(w)+\hat q_{*w}(w)\in{\frak a}$. It is clear that
$d\Phi=\frac12id\tilde\theta$ if
the first term in the right-hand side of~(4.35)
equals $\frac12i\tilde\theta$ and the second one determines a closed
1-form $\sigma$ on $W_{\bold a}^\Sigma$. In other words, if
$(1+SR^{-1}SR^{-1})_wR_w(w')=w$,
i.e.
$$
b(x)(2xq(x^2)+2x^3q'(x^2))=b(x)\bigl(x^2q(x^2)\bigr)'=x
$$
because of relation~(2.4) and~(4.9). Next, the form
$$
\sigma_w(u)=2^{-1}\langle (1+iSR^{-1}_w)
(\hat q(w)+\hat q_{*w}(w)),u\rangle
$$
is closed. Indeed,  $SR^{-1}_w=a_1\Upsilon_{*w}$. The mapping
$w\mapsto (1+ia_1\Upsilon_{*w})(\hat q(w)+\hat q_{*w}(w))$ is
$K$-equivariant because so are the mappings
$\hat q$,
$w\mapsto \hat q_{*w}$ and $w\mapsto \Upsilon_{*w}$.
Now by Lemma~4.5, $d\sigma=0$ because a restriction
$\sigma$ to $W\cap{\frak a}$ is a linear combination of
the closed forms
$(1-i/{x_j^2})d\bigl(x_j^2q(x_j^2)\bigr)$.
{}\hfill$\square$
\enddemo
\example{Example 4.14}
{\it
Hyperk\"ahler structures on homogeneous domains
in $T(SU(2)/U(1))$}.
Let $G/K=SU(2)/U(1)$. The Lie algebra
${\frak g}=su(2)$ has the basis $X_1,Y_1,T_1$~(4.8).
For each complex number
$z=x+\imath y$ put
$zX_1\overset\text{def}\to{=} xX_1+yIX_1=xX_1+yY_1\in{\frak m}$.
Let
$(J(R+iS),\Omega)$ be a K\"ahler structure
anticommuting with $J^-$.
Then
for each $w=zX_1$, $|z|^2>a_\dagger$
$$
R_{w}=\psi(|z|)\cdot\operatorname{Id}_{\frak m},
\ \text{where}\qquad
\psi(x)=\frac{\sqrt{x^6\cdot(x^4+a_0x^2-(a_1^2+a_2^2))}}
{x^4+(a_1^2+a_2^2)},
$$
and
$$
S_{w}(vX_1)=\psi(|z|)(a_1+\imath a_2)\overline{z^{-2}v}\cdot X_1
\qquad\text{for all}\qquad v\in{\Bbb C}.
$$
\endexample

\remark{Remark 4.15}
The global hyper-K\"ahlerian structure
$({\bold g},J^-, J(P))$, where
$J(P)$ is defined by~(4.25) with
$a_0>0,\varepsilon=1$ and
$a_1=a_2=0$, coincides with the structure
constructed in~\cite{BG1, Theorem 1}. Here
${\bold g}$ is the corresponding to the pair
$(J(P),\Omega)$ (hyper)K\"ahler metric. This
is a unique metric for which its restriction to
$G/K$ coincides with the metric
$\sqrt{a_0}\cdot{\bold g}_{\text M}$ on
$M=G/K$ (because
$P_0=\sqrt{a_0}\cdot\operatorname{Id}_{\frak m}$).
Our formula for a potential generalizes
such a formula in~\cite{BG2}.
\endremark

\enddocument